\documentclass{siamltex}
\usepackage{amssymb}
\usepackage{bm} 
\usepackage{graphicx}
\usepackage{graphics}
\usepackage{caption2}
\usepackage{psfrag}
\usepackage{amsmath}
\usepackage{amscd}
\usepackage{amsfonts}
\usepackage{float}
\usepackage{latexsym}
\usepackage{lscape}
\usepackage{color}
\usepackage{lineno}
\usepackage{algorithmicx,algorithm}

\newtheorem{remark}{Remark}[section]
\newtheorem{assumption}{Assumption}[section]

\newcommand{\ds}{\displaystyle}

\begin{document}

\title{Asymptotic stability of high-index saddle dynamics \thanks{The work was partially supported by National Natural Science Foundation of China Nos. 12225102, T2321001, 12288101, and 12426653 (for LZ); the National Natural Science Foundation of China Nos. 12301555 and 12622126, and Natural Science Foundation of Shandong Province No. ZR2025QB01 (for XZ).}}
\author{
Baoming Shi\thanks{Department of Applied Physics and Applied Mathematics,  Columbia University, New York, NY 10027, USA. 
(bs3705@columbia.edu)}
  \and Lei Zhang\thanks{School of Mathematical Sciences, Beijing International Center for Mathematical Research, Center for Quantitative Biology, Center for Machine Learning Research, Peking University, Beijing, 100871, China. (zhangl@math.pku.edu.cn)}
  \and Xiangcheng Zheng\thanks{School of Mathematics, Shandong University, 
    Jinan, 250100, China.     (xzheng@sdu.edu.cn; sqzhu@mail.sdu.edu.cn)}
  \and Shangqin Zhu\footnotemark[4]
}
\maketitle

\begin{abstract}
High-index saddle dynamics (HiSD) has been extensively adopted to compute saddle points  of arbitrary index, with various applications such as constructing solution landscapes of complex systems. In this work, we prove the asymptotic stability of various HiSD methods, including the idealized, perturbed, and coupled HiSD methods, as well as their momentum-accelerated variants. The analysis is Jacobian-free  and thus compensates for the linear stability analysis in the literature.  The results not only provide specific characterizations for the convergence rate and the basin of attraction, but also reveal intrinsic correlations among dynamics in the coupled HiSD. Theoretical results are substantiated by numerical experiments.
\end{abstract}



\pagestyle{myheadings}
\thispagestyle{plain}

\markboth{B. Shi, L. Zhang, X. Zheng and S. Zhu}{Asymptotic stability of high-index saddle dynamics}

\section{Introduction}
Saddle points play a critical role in various fields. In particular, the energy and solution landscapes of complex systems, which contain saddle points and local minimums and the connections among them, have been successfully applied in studying, e.g., defect configurations in liquid crystals \cite{wang2021modeling},  cell fate decision \cite{ZhaLi}, gene expression and metabolism \cite{Chaves}, nucleation of quasicrystals \cite{YinPNAS}, and hierarchical and ultrametric barriers in jammed granular matter \cite{Wu}.  The saddle points can be classified by the Morse index, which is characterized by the maximal dimension of a subspace on which the Hessian is negative definite \cite{Milnor}.
The computation of saddle points is challenging due to their unstable nature, and there exist several successful numerical algorithms, see e.g. \cite{Bao,Che,Chenliu,Doye,EZho,Farr, Farr2,Gao,Gou,Li2001,SanWan,ZhaDu}.

In \cite{Bof,Chu,YinSISC}, a coupled high-index saddle dynamics (HiSD) approach is proposed to find any-index saddle points, whose mechanism lies in searching along a combined direction according to the manifold property of the target saddle point. The linear stability of the coupled HiSD was analyzed in \cite{Chu,YinSISC}, following a similar analysis for the index-1 dimer dynamics presented earlier in \cite{ZhaDu}. It has since then been extended to an accelerated HiSD \cite{Luoacc} and other variants \cite{du2026convergence}. For the asymptotic stability without linearization,  the discrete versions of the idealized HiSD (i.e., the unstable subspace is assumed to be exact at each time instant) and the perturbed HiSD by inexact unstable subspace are analyzed in \cite{Luo}. Later, the asymptotic stability is also considered for the discrete version of an accelerated HiSD \cite{Luoacc}. 

Despite the above progress on the stability analysis of HiSD, the asymptotic stability analysis of continuous formulations of HiSD, including the idealized, perturbed, and coupled HiSD methods, as well as their accelerated variants, remains untreated. 
Although the known linear stability results characterize the  stability of the nonlinear dynamics locally, it is difficult to give a specific quantification for the basin of attraction. 
This motivates the current work on the original asymptotic stability of HiSD without linearization. The contributions and novelties are summarized as follows:
\begin{itemize}
 \item For the analysis of idealized HiSD, we combine spectral analysis with the energy method. Under local Lipschitz assumptions, we establish the boundedness of the dynamics through a ``local-to-global'' procedure and then apply a post-processing step to obtain sharper estimates (cf. Theorem \ref{T1}). For the perturbed HiSD, since the inexact unstable subspace is no longer spanned by eigenvectors, the spectral properties of the Hessian within this subspace are not apparent. To address this challenge, we employ a matrix splitting strategy that decomposes the Hessian into its dominant spectral structure and perturbation terms arising from subspace inaccuracy. This allows us to quantify how unstable subspace inaccuracy affects the basin of attraction and the convergence rate (cf. Theorem \ref{i1}).
 
 \item For the coupled HiSD, the main difficulty lies in the simultaneous evolution of the state variable and the unstable subspace such that no a priori estimate of the unstable subspace is available. Furthermore, unlike the linearized formulation studied in \cite{Chu,YinSISC}, the dynamics of the unstable subspace does not exhibit favorable structures for analysis. Thus, we first derive a priori bound on the unstable subspace projection error in the Frobenius norm (cf.  Lemma \ref{L8}). Then we combine this with the state dynamics to establish asymptotic stability by jointly controlling the state and unstable-subspace errors, ensuring that neither exceeds its initial bound throughout the evolution. An interesting phenomenon is that the convergence rate of the state variable does not depend on that of the unstable subspace, otherwise the opposite (cf. Theorem \ref{i3}). These reveal intrinsic correlations among dynamics in the coupled HiSD, as validated by numerical experiments (cf. Section \ref{sec62}).
 
 \item For the momentum-accelerated variants of idealized, perturbed and coupled HiSD, we reformulate the second-order flows as first-order systems. Since the resulting coefficient matrices are asymmetric, we adopt an integrating factor approach \cite{Wal} and derive matrix exponential bounds through spectral analysis. Combining these bounds with estimates of the nonlinear remainder terms, we establish asymptotic stability through a continuation argument, obtaining explicit bounds on the basin of attraction and the convergence rate (cf. Theorems \ref{aa}--\ref{i4}). Based on the theoretical results, we derive explicit sufficient conditions for faster convergence and numerically examine the feasible parameter regions (cf. Section \ref{sec61}). These results provide practical guidance for achieving acceleration.
 \end{itemize} 
 
The rest of this paper is organized as follows: In Section \ref{sec2} we present formulations of different HiSD methods. In Sections \ref{idea}--\ref{acce}, the idealized and perturbed HiSD, the coupled HiSD, and their accelerated variants are analyzed, respectively. Section \ref{num} presents numerical experiments, and the last section concludes the paper.

\section{Formulations of HiSD}\label{sec2}
Given a twice Fr\'echet differentiable energy functional $E(x)$ with $x\in\mathbb R^N$, the idealized HiSD for finding an index-$k$ saddle point $x^*$ of $E(x)$ reads \cite{Chu,EZho,Qua,YinSISC}
\begin{equation}\label{sadk}
\begin{array}{l}
\ds \frac{dx}{dt} =-\beta(I-2\hat V\hat V^\top) \nabla E(x).
\end{array}
\end{equation}
Here $x=x(t)$ is the state variable, $\nabla E$ refers to the gradient of $E$, $\hat V:=[\hat v_1,\ldots,\hat v_k]\in\mathbb R^{N\times k}$ where $\{\hat v_1,\ldots,\hat v_k\}$ are orthonormal eigenvectors of the Hessian $\nabla^2 E(x)$ corresponding to its first $k$ smallest eigenvalues, and $\beta>0$ is a relaxation parameter. In other words, the $\hat V$ corresponds to the information of the unstable subspace such that, acted by the Householder transform, the right-hand side of (\ref{sadk}) is the combination of the ascent direction on the unstable subspace and the descent direction on its stable  orthogonal complement space. Such searching mechanism is consistent with the local manifold property of the target saddle point, and thus ensures the effectiveness of (\ref{sadk}) in locating saddle points.

In the idealized HiSD (\ref{sadk}), it is assumed that $\{\hat v_i\}_{i=1}^k$ are exact orthonormal eigenvectors of the Hessian $\nabla^2 E(x)$ corresponding to its $k$ smallest eigenvalues. In practice, $\hat V$ may contain certain errors due to, e.g., numerical computation, leading to inexact eigenvectors collected in $\tilde V$. Then the corresponding perturbed HiSD take the following form
\begin{equation}\label{iesd}
\ds \frac{dx}{dt} =-\beta  (I-2\tilde V\tilde V^\top)  \nabla E(x).
\end{equation}

In the above HiSD formulations, evaluating eigenvectors at each time step can be computationally expensive, especially for large-scale problems. An alternative approach is to treat the eigenvectors as dynamic variables that evolve toward the true eigenvectors by optimizing the Rayleigh quotient. This leads to the coupled HiSD where both the state variable $x$ and the unstable subspace $V:=[v_1,\ldots,v_k]\in\mathbb R^{N\times k}$ with $V^\top V=I$ are updated simultaneously \cite{Bof,Chu,YinSISC}:
\begin{equation}\label{hisd}
    \left\{
    \begin{aligned}
    \frac{dx}{dt} & = -\beta(I-2VV^\top) \nabla E(x),\\[0.05in]
    \frac{dV}{dt} &= -\gamma\big(I-VV^\top\big) \nabla^2 E(x)V.
    \end{aligned}
    \right.
\end{equation}
Here $\gamma >0$ is a  relaxation parameter. It is proved in \cite{Chu,YinSISC} that a linearly stable stationary point of (\ref{hisd}) is an index-$k$ saddle point.

Recently, an accelerated HiSD was proposed in \cite{Luoacc} to improve the convergence behavior by introducing an additional momentum term, which results in the following second-order flow  
 \begin{equation}\label{asadk}
 \left\{
    \begin{aligned}
   \frac{dx}{dt} &=m,\\[0.05in]
    \frac{dm}{dt}&=-\alpha_1 m-\alpha_2(I-2\hat V\hat V^\top) \nabla E(x),
    \end{aligned}
\right.
\end{equation}
 where $\alpha_1,\alpha_2>0$ are relaxation parameters. Correspondingly, we consider the perturbed HiSD with momentum acceleration
\begin{equation}\label{ieasadk}
 \left\{
    \begin{aligned}
   \frac{dx}{dt} &=m,\\[0.05in]
    \frac{dm}{dt}&=-\alpha_1 m-\alpha_2(I-2\tilde V\tilde V^\top) \nabla E(x),
    \end{aligned}
\right.
\end{equation}
and the coupled HiSD with momentum acceleration
\begin{equation}\label{ieah}
 \left\{
    \begin{aligned}
   \frac{dx}{dt} &=m,\\[0.05in]
    \frac{dm}{dt}&=-\alpha_1 m-\alpha_2(I-2 V V^\top) \nabla E(x),\\[0.05in]
    \frac{dV}{dt} &= -\gamma\big(I-VV^\top\big) \nabla^2 E(x)V.
    \end{aligned}
\right.
\end{equation}


Let $\|\cdot\|$ and $\|\cdot\|_F$ be the standard $l^2$ and Frobenius norms of a vector or a matrix, respectively. Now we introduce the concept of exponential asymptotic stability \cite{Wal}.

\begin{definition} (Exponential asymptotic stability) A constant solution $x^*$ of a first-order system of $x$ is said to be \textit{exponential asymptotic stable} if there exist some $K,\sigma,\delta_2>0$ such that $\|x(0)-x^*\|< \delta_2$ implies
 $$\|x(t)-x^*\|\leq Ke^{-\sigma t},~~t\geq 0.$$
\end{definition}

Throughout this work, we impose the following standard conditions \cite{Luo,Luoacc}:

\begin{assumption}
The solutions to the proposed HiSD are continuous over $t\in[0,\infty)$, and there exists a neighborhood $B(x^*,\rho) =\{x|\|x-x^*\|\leq\rho\}$ of the target saddle point $x^*$ for some $\rho>0$ such that:
    \begin{enumerate}
    \item[(i)] The Lipschitz condition
     $\|\nabla^2E(x)-\nabla^2E(y)\|\leq M\|x-y\|$
     holds for $ x,y\in B(x^*,\rho)$ and for some $M>0$;
   \item[(ii)] For any $x\in B(x^*,\rho)$, the eigenvalues $\{\lambda_i\}_{i=1}^N$ of $\nabla^2E(x)$ satisfy
   $\lambda_1\leq\ldots\leq \lambda_k<0< \lambda_{k+1}\leq\ldots\leq\lambda_N$
         and $|\lambda_i|\in[\mu,L]$ for $1\leq i\leq N$ for some positive constants $0<\mu< L$.
\end{enumerate}
\label{ass1}
\end{assumption}

\section{Stability of idealized  and perturbed HiSD}\label{idea}
In this section, we establish the asymptotic stability of both idealized and perturbed HiSD.
\subsection{Analysis of idealized HiSD} 
We apply $\nabla E(x^*)=0$ to rewrite (\ref{sadk}) as
 \begin{align}
\frac{dx}{dt}& =-\beta(I-2\hat V\hat V^\top) (\nabla E(x)-\nabla E(x^*)) \nonumber\\
&=-\beta (I-2\hat V\hat V^\top) \int_0^1\nabla^2 E(x^*+s(x-x^*))ds(x-x^*)\label{11}\\
&=- \beta (I-2\hat V\hat V^\top)\nabla^2 E(x)(x-x^*)+q(x)\nonumber
\end{align}
 where
 $$q(x):=-\beta(I-2\hat V\hat V^\top) \int_0^1\big[\nabla^2 E(x^*+s(x-x^*))-\nabla^2 E(x)\big]ds(x-x^*).$$

\begin{theorem}\label{T1}
Suppose Assumption \ref{ass1} holds and
\begin{equation}\label{x0}
\|x(0)-x^*\|< \min\Big\{\rho,\frac{2\mu}{M}\Big\}.
\end{equation}
 Then the target index-$k$ saddle point $x^*$ of the idealized HiSD (\ref{sadk}) is exponential asymptotic stable. In particular, for $\varepsilon>0$ satisfying $\varepsilon<\frac{\|x(0)-x^*\|M}{2}, $
 the following estimate holds for some $K(\varepsilon)>0$
\begin{align*}
\|x-x^*\|\leq K(\varepsilon) e^{-\beta(\mu-\varepsilon) t},~~t\geq 0.
\end{align*}
\end{theorem}
\begin{proof}
The proof will be divided into two parts.

\textbf{Part I: Boundedness of $x$.}
Due to (\ref{x0}) and the continuity of the solutions, there exists a $\delta_0>0$ such that
$$x(t)\in B_0:=\Big\{x:\|x-x^*\|\leq  \min\Big \{\rho,\displaystyle\frac{2\mu}{M} \Big\}\Big\},~~t\in [0,\delta_0].$$
Thus, for $t\in [0,\delta_0]$, we apply the Assumption \ref{ass1} and $\|(I-2\hat V\hat V^\top)v\|=\|v\|$ for any $v\in\mathbb R^N$ to get
 \begin{equation*}
 	\|q\|\leq \beta M\|x-x^*\|^2\int_0^1(1-s)ds=\frac{\beta M}{2}\| x-x^*\|^2.
 \end{equation*}
 By eigenvalue decomposition theorem and the Assumption \ref{ass1}, we have 
 $\nabla^2 E(x)=\sum_{i=1}^N\lambda_i\hat v_i \hat v_i^\top,$
 which implies
 \begin{equation*}
 \begin{aligned}
-(I-2\hat V\hat V^\top)\nabla^2 E( x)&=\Big(I -2\sum_{j=1}^k \hat v_j \hat v_j^\top\Big)\Big(-\sum_{i=1}^N\lambda_i \hat v_i \hat v_i^\top\Big)\nonumber\\
& =\sum_{j=1}^k\lambda_j \hat v_j \hat v_j^\top-\sum_{j=k+1}^N\lambda_j\hat v_j \hat v_j^\top= \mathcal{V}  \Lambda \mathcal{V}^\top, 
 \end{aligned}	
 \end{equation*}
 where $\mathcal{V}=[\hat v_1,\ldots,\hat v_N]$ is an orthogonal matrix and $\Lambda$ is a diagonal matrix with entries $\Lambda_i$ for $1\leq i\leq N$ defined by
 $$\Lambda_i=
 \left\{
 \begin{array}{ll}
 \lambda_i,~~&1\leq i\leq k;\\
   -\lambda_{i},~~&k+1\leq i\leq N.
 \end{array}
 \right.
 $$
 By Assumption \ref{ass1}, we have $\Lambda_i\leq -\mu$ for $1\leq i\leq N$. Then we multiply both sides of (\ref{11}) by $ x-x^*$, and invoke the above estimates to get
 \begin{align}
 &\frac{1}{2}\frac{d\|x-x^*\|^2}{dt}=(x-x^*)^\top\frac{d x}{dt}=\beta (x-x^*)^\top\mathcal{V}\Lambda \mathcal{V}^\top (x-x^*)+(x-x^*)^\top q\nonumber \\
&\qquad= \beta\sum_{i=1}^N\Lambda_i[\mathcal{V}^\top(x-x^*)]_i^2+(x-x^*)^\top q\label{dx}\\
&\qquad\leq -\beta\mu\|\mathcal{V}^\top(x-x^*)\|^2+\frac{\beta M}{2}\|x-x^*\|^3 =-\beta(\mu-\frac{ M}{2}\|x-x^*\|)\|x-x^*\|^2,\nonumber
 	\end{align}
where $[\mathcal{V}^\top(x-x^*)]_i$ represents the $i$th entry of $\mathcal{V}^\top(x-x^*)$ and we used $\|\mathcal{V}^\top(x-x^*)\|=\|x-x^*\|$ since $\mathcal{V}$ is an orthogonal matrix. We cancel $\|x-x^*\|$ on both sides to get
\begin{equation}\label{ls2}
\frac{d\|x-x^*\|}{dt}\leq -\beta\big(\mu-\frac{M}{2}\|x-x^*\|\big)\|x-x^*\|. 
\end{equation}
By (\ref{x0}) we have
\begin{equation}\label{14}
\|x(0)-x^*\|< \frac{2\mu}{M}\Longrightarrow \mu-\frac{M\|x(0)-x^*\|}{2}>0.
\end{equation}
Then (\ref{ls2}) implies
$$\frac{d\|x-x^*\|}{dt}\Big |_{t=0^+}\leq -\beta\big(\mu-\frac{M}{2}\|x(0)-x^*\|\big)\|x(0)-x^*\|<0. $$ 
Thus, there exists a $0<\delta\leq \delta_0$ such that $\|x-x^*\|<\|x(0)-x^*\|$ for any $0< t\leq \delta$. 

Now we consider the set
$U=\{t\geq 0:\|x-x^*\|>\|x(0)-x^*\|\}$. We intend to prove that $U$ is empty. If not, then $U$ is a non-empty set with a lower bound $0$. Thus $U$ has an infimum denoted by $t_0$. Since we have proved that $\|x-x^*\|<\|x(0)-x^*\|$ for any $0< t\leq \delta$, $t_0$ can not be less than $\delta$ and is thus positive.  As $\|x-x^*\|\leq \|x(0)-x^*\|$ for $t\in [0,t_0]$, the (\ref{ls2}) holds for $t\in [0,t_0]$ such that we apply (\ref{ls2}) and (\ref{14}) to obtain
\begin{align*}
&\frac{d\|x-x^*\|}{dt}\Big|_{t=t_0}\leq -\beta\big(\mu-\frac{M}{2}\|x(t_0)-x^*\|\big)\|x(t_0)-x^*\|\\
&\qquad= -\beta\big(\mu-\frac{M}{2}\|x(0)-x^*\|\big)\|x(0)-x^*\|<0.
\end{align*}
This leads to $\|x(t_0-t)-x^*\|>\|x(t_0)-x^*\|=\|x(0)-x^*\|$ for $0<t\leq \epsilon $ for some small enough $\epsilon$, which, according to the definition of $t_0$, implies $t_0\leq t_0-\epsilon$, a contradiction to $\epsilon>0$. 

 Based on the above derivation, we conclude that $U$ is empty such that $\|x-x^*\|\leq \|x(0)-x^*\|$ for any $t\geq 0$. We combine this with (\ref{x0}) to conclude that $ x(t)\in B_0$ for any $t\geq 0$ such that (\ref{ls2}) holds for any $t\geq 0$.

\textbf{Part II: Stability analysis.}
We rewrite (\ref{ls2}), which is valid for all $t>0$ due to $x\in B_0$, by the integrating factor method
 $$\frac{d}{dt}\Big(\exp\big(\int_0^t \beta(\mu-\frac{M}{2}\|x(s)-x^*\|) ds\big)\|x-x^*\|\Big)\leq 0.$$
 Integrating this equation from $0$ to $t$ yields
\begin{align}\label{ls3}
\|x-x^*\| \leq \|x(0)-x^*\|\exp\Big(-\int_0^t\beta(\mu-\frac{M}{2}\|x(s)-x^*\|)ds\Big) .
\end{align}
We first apply $\|x-x^*\|\leq \|x(0)-x^*\|$ to obtain a preliminary result from (\ref{ls3}) 
\begin{align*}
\|x-x^*\|& \leq \|x(0)-x^*\|\exp\Big(-\int_0^t\beta(\mu-\frac{M}{2}\|x(0)-x^*\|)ds\Big) \\
&= \|x(0)-x^*\|\exp\Big(-\beta(\mu-\frac{M}{2}\|x(0)-x^*\|)t\Big) .
\end{align*}
Then for $\varepsilon>0$ satisfying
$\varepsilon<\frac{\|x(0)-x^*\|M}{2}, $
we have 
\begin{equation}\label{ls5}	
\|x-x^*\| \leq\|x(0)-x^*\|\exp\Big(-\beta(\mu-\frac{M}{2}\|x(0)-x^*\|)t\Big) \leq  \frac{2\varepsilon}{M}, 
\end{equation}
\vspace{-0.2in}
\begin{equation}\label{ls6}
t\geq t_\varepsilon:=\frac{\ds -\ln \frac{2\varepsilon}{\|x(0)-x^*\|M}}{\ds \beta(\mu-\frac{M}{2}\|x(0)-x^*\|)}. 
\end{equation}
We then invoke (\ref{ls5})--(\ref{ls6}) in (\ref{ls3}) to get a sharper decay rate for $t\geq t_{\varepsilon}$
\vspace{-0.05in}
\begin{align*}
\|x-x^*\| &\leq \|x(0)-x^*\|\exp\Big(-\Big[\int_0^{t_\varepsilon}+\int_{t_\varepsilon}^t\Big]\beta(\mu-\frac{M}{2}\|x(s)-x^*\|)ds\Big)\\
&\leq \|x(0)-x^*\|\exp\Big(-\beta(\mu-\frac{M}{2}\|x(0)-x^*\|)t_\varepsilon\Big) \exp\Big(-\int_{t_\varepsilon}^t\beta(\mu-\varepsilon)ds\Big) \\
&= \frac{2\varepsilon}{M}e^{\beta(\mu-\varepsilon)t_{\varepsilon}}e^{-\beta(\mu-\varepsilon)t}.
\end{align*}
For $0\leq t\leq t_\varepsilon$, we could apply a rough bound as
\begin{equation*}\label{ls8}
\begin{aligned}
 \|x-x^*\|&=e^{\beta(\mu-\varepsilon)t}\|x-x^*\|e^{-\beta(\mu-\varepsilon)t} \leq e^{\beta(\mu-\varepsilon)t_\varepsilon}\|x(0)-x^*\|e^{-\beta(\mu-\varepsilon)t}.
 \end{aligned} 
\end{equation*}
We combine the above estimates to complete the proof.
\end{proof}

\subsection{Analysis of perturbed HiSD}
To analyze the perturbed HiSD by inexact unstable subspace, we require a proper distance measure between the exact subspace spanned by $\hat V$ and the approximate subspace spanned by $\tilde{V}$. We evaluate this distance from a projection perspective and assume that there exists an $\alpha \in [0,1]$ such that 
$\|\hat V {\hat V}^\top - \tilde{V} \tilde{V}^\top \|_2 \leq  \alpha.$
Similar to (\ref{11}), we apply $\nabla E(x^*)=0$ to rewrite the perturbed HiSD (\ref{iesd}) as
\begin{equation*}
\frac{dx}{dt}= -\beta (I-2\hat V \hat V^\top) \nabla^2 E(x)(x-x^*)+ r(x)+\tilde q(x) ,
\end{equation*}
where
\begin{equation*}
	\begin{gathered}
r(x):= -\beta \big((I-2\tilde V \tilde V^\top)-(I-2\hat V \hat V^\top) \big)\nabla^2 E(x)(x-x^*),\\
\tilde q(x):=-\beta (I-2\tilde V\tilde V^\top)  \int_0^1\big[\nabla^2 E(x^*+s( x-x^*))-\nabla^2 E(x)\big]ds(x-x^*).
	\end{gathered}
\end{equation*}
Compared with (\ref{11}), the only difference lies in the additional term $r$, the $l^2$ norm of which can be bounded by $2\alpha \beta L\|x-x^*\|$ if $x\in B(x^*,\rho)$. Then a similar derivation as (\ref{dx}) holds for $x\in B(x^*,\rho)$
\begin{equation*}
 \begin{aligned}
 \frac{1}{2}\frac{d\|x-x^*\|^2}{dt}&=(x-x^*)^\top\frac{d x}{dt}=\beta (x-x^*)^\top\mathcal{V}\Lambda \mathcal{V}^\top (x-x^*)+(x-x^*)^\top (r+ \tilde q)\\
&= \beta\sum_{i=1}^N\Lambda_i[\mathcal{V}^\top(x-x^*)]_i^2+(x-x^*)^\top r+(x-x^*)^\top \tilde q\\
&\leq -\beta\mu\|\mathcal{V}^\top(x-x^*)\|^2+ 2\alpha \beta L\|x-x^*\|^2 +\frac{\beta M}{2}\|x-x^*\|^3 \\
& =-\beta\Big(\hat\mu-\frac{ M}{2}\|x-x^*\|\Big)\|x-x^*\|^2,
 	\end{aligned}
 \end{equation*}
where $\hat \mu:=\mu-2\alpha L$ and $\alpha$ is selected such that $\hat\mu>0$.  Thus, with $\mu$ replaced by $\hat\mu$ in the proof of Theorem \ref{T1}, one can prove the exponential asymptotic stability of the perturbed HiSD (\ref{iesd}) (we omit the details due to similarity).

\begin{theorem}\label{i1} Suppose Assumption \ref{ass1} holds, and $\alpha$ is small enough such that $\hat \mu:=\mu-2\alpha L>0$. Then for $0<\varepsilon<\hat \mu$ and 
$
 	\|x(0) - x^* \| < \min \big\{\rho, \frac{2 \hat \mu}{M} \big\}, 
$
the following estimate holds for perturbed HiSD (\ref{iesd}) for some $\tilde K(\varepsilon)>0$
\begin{equation*}
\|x-x^*\|\leq  \tilde K(\varepsilon) e^{-\beta(\hat \mu -\varepsilon) t},~~t\geq 0.
\end{equation*}
\end{theorem}
\section{Stability of coupled HiSD}\label{coup} We then consider the  asymptotic stability of the coupled HiSD (\ref{hisd}). Compared with the above uncoupled cases, the main difficulty of the coupled case lies in the lack of a priori estimate and favorable dynamical structures of the unstable subspace. Thus, we first derive a priori bound for the asymptotic stability of the projection operator to the unstable subspace, and then combine with the dynamics of the state variable  to perform the analysis. In this analysis, we use the Frobenius norm to facilitate differentiation and employ proofs by contradiction to address the coupling between the dynamics.

\subsection{Auxiliary estimates} We first give a priori estimate for the projection $VV^\top$.
\begin{lemma}\label{L8}
Suppose the Assumption~\ref{ass1} holds. If the solutions $\{x, V\}$ to \eqref{hisd} satisfy $x(t) \in B(x^*, \rho)$ and $\|VV^\top - V^*{V^*}^\top\|_F \leq 1$ at some time instant, then the following inequality holds at this instant
\begin{equation*}
    \frac{d}{dt}\|VV^\top - V^*{V^*}^\top\|_F
    \leq -\gamma\mu\|VV^\top - V^*{V^*}^\top\|_F
    + 2\gamma \sqrt{k} M\|x - x^*\|.
\end{equation*}
\end{lemma}

\begin{proof}
Direct computations give
\begin{equation}\label{half}
\begin{aligned}
\frac{d}{dt}\|VV^\top - V^*{V^*}^\top\|_F^2
&= 2\bigl\langle VV^\top - V^*{V^*}^\top,\,
   \dot{V}V^\top + V\dot{V}^\top \bigr\rangle_F\\
&= 4\bigl\langle (VV^\top - V^*{V^*}^\top)V,\, \dot{V} \bigr\rangle_F.
\end{aligned}
\end{equation}
We invoke $\dot{V} = -\gamma(I - VV^\top)\nabla^2 E(x)V$ and decompose the resulting equation to get
\begin{equation*}
\begin{aligned}
&4\bigl\langle (VV^\top - V^*{V^*}^\top)V,\, \dot{V} \bigr\rangle_F \\
&\qquad =4\bigl\langle (VV^\top - V^*{V^*}^\top)V,\;
   \gamma(I-VV^\top)(\nabla^2 E(x^*)-\nabla^2 E(x))V\bigr\rangle_F\\
&\qquad \quad + 4\bigl\langle (VV^\top - V^*{V^*}^\top)V,\;
   -\gamma(I-VV^\top)\nabla^2 E(x^*)V\bigr\rangle_F =: A_1 + A_2.
\end{aligned}
\end{equation*}

By  $\|\nabla^2 E(x) - \nabla^2 E(x^*)\| \leq M\|x - x^*\|$, $\|(I - VV^\top)W\|_F \leq \|W\|_F$ for any $W\in\mathbb R^{N\times k}$,  $\|\Sigma  V\|_F \leq \|\Sigma\|\|V\|_F$ for any $\Sigma \in\mathbb R^{N\times N}$ and $\|V\|_F = \sqrt{k}$, we bound $A_1$ as
\begin{equation*}
\begin{aligned}
|A_1|
&\leq 4\gamma\|(VV^\top - V^*{V^*}^\top)V\|_F
       \cdot\|(I - VV^\top)(\nabla^2 E(x) - \nabla^2 E(x^*))V\|_F \\
&\leq 4\gamma\|(VV^\top - V^*{V^*}^\top)V\|_F
       \cdot\|\nabla^2 E(x) - \nabla^2 E(x^*)\|\cdot\|V\|_F \\
&\leq 4 \gamma \sqrt{k} M\|(VV^\top - V^*{V^*}^\top)V\|_F\cdot\|x - x^*\|.
\end{aligned}
\end{equation*}
Using $\|(VV^\top - V^*{V^*}^\top)V\|_F \leq \|VV^\top - V^*{V^*}^\top\|_F$ (since $\|V\|_2 = 1$) yields
\begin{equation}\label{A1}
|A_1| \leq 4 \gamma \sqrt{k}M\|VV^\top - V^*{V^*}^\top\|_F\cdot\|x - x^*\|.
\end{equation}

To bound $A_2$, set $S := {V^*}^\top V$ and let $V^*_\perp = [v_{k+1}^*, \dots, v_N^*]$
be the orthogonal complement of $V^*$. Then
$
V = V^* S + V^*_\perp Y$ with $ Y := {V^*_\perp}^\top V$
and $S^\top S + Y^\top Y = I$. Using the eigenrelations $\nabla^2 E(x^*)V^* = V^*\Lambda^*$ and $ \nabla^2 E(x^*)V^*_\perp = V^*_\perp \Lambda^*_\perp$, where
$\Lambda^* = \diag(\lambda_1^*, \dots, \lambda_k^*)$ and $\Lambda^*_\perp = \diag(\lambda_{k+1}^*, \dots, \lambda_N^*)$, we obtain
\begin{equation*}
V^\top\nabla^2 E(x^*)V =  (V^* S + V^*_\perp Y)^\top\nabla^2 E(x^*)(V^* S + V^*_\perp Y) = S^\top\Lambda^* S + Y^\top\Lambda^*_\perp Y.
\end{equation*}
We apply this to obtain
\begin{equation}\label{A2e}
\begin{aligned}
A_2& = 4\bigl\langle (VV^\top - V^*{V^*}^\top)V,\;-\gamma(I-VV^\top)\nabla^2 E(x^*)V\bigr\rangle_F\\
& =4 \big\langle V - V^* S , -\gamma(I- VV^\top) \nabla ^2 E (x^*)V \big\rangle_F\\
& = 4 \big\langle V^* S , \gamma(I- VV^\top) \nabla ^2 E (x^*)V \big\rangle_F\\
&=4\gamma\operatorname{tr}\bigl(S^\top {V^*}^\top {\nabla^2 E(x^*)}^\top V
   - S^\top {V^*}^\top V (V^\top\nabla^2 E(x^*)V)\bigr) \\
&=4\gamma\operatorname{tr} \bigl(S^\top (V^*\Lambda^*)^\top V
   - S^\top S (V^\top\nabla^2 E(x^*)V)\bigr) \\
& =4\gamma\operatorname{tr} \bigl(S^\top\Lambda^* S
   - S^\top S(S^\top\Lambda^* S + Y^\top\Lambda^*_\perp Y)\bigr)\\
&= 4\gamma\operatorname{tr} \bigl((I - S^\top S)S^\top\Lambda^* S\bigr)
  -4\gamma\operatorname{tr} \bigl(S^\top S Y^\top\Lambda^*_\perp Y\bigr).
\end{aligned}
\end{equation}
Under Assumption~\ref{ass1}, $\lambda_i^* \leq -\mu < 0$ for $1 \leq i \leq k$. Thus, both $ I - S^\top S $ and $S^\top (-\mu I-\Lambda^*)S$ are symmetric and positive semidefinite such that
\begin{equation}\label{trace1}
\operatorname{tr} ((I - S^\top S)S^\top (-\mu I-\Lambda^*)S)\geq 0 \Rightarrow \operatorname{tr} \bigl((I - S^\top S)S^\top\Lambda^* S\bigr) \leq -\mu\operatorname{tr} \bigl((I - S^\top S)S^\top S\bigr).
\end{equation}
Similarly, since $\lambda_j^* \geq \mu$ for $j > k$, both $ S^\top S $ and $Y^\top (\Lambda^*_\perp -\mu I)Y$ are symmetric and positive semidefinite such that
\begin{equation}\label{trace2}
\operatorname{tr} \bigl(S^\top S Y^\top (\Lambda^*_\perp-\mu I) Y\bigr)
\geq 0 \Rightarrow 
\operatorname{tr} \bigl(S^\top S Y^\top\Lambda^*_\perp Y\bigr)
\geq \mu\operatorname{tr} \bigl(S^\top S Y^\top Y\bigr).
\end{equation}
Substituting \eqref{trace1} and \eqref{trace2} into \eqref{A2e} and using
$S^\top S + Y^\top Y = I$, we arrive at 
\begin{equation*}
A_2 \leq -4\gamma\mu\operatorname{tr} \bigl((I - S^\top S)S^\top S\bigr)
         -4\gamma\mu\operatorname{tr} \bigl(S^\top S Y^\top Y\bigr)
    = -8\gamma\mu\operatorname{tr}(S^\top S Y^\top Y). 
\end{equation*}
The trace identity $\operatorname{tr}(P(M))=\sum_{i=1}^k P(\lambda_i(M))$
for any polynomial $P$ gives $A_2 \leq -8\gamma\mu\sum_{i=1}^k\lambda_i^S\lambda_i^Y$, where $\{\lambda_i^{S}\}_{i=1}^k$ and $\{\lambda_i^{Y}\}_{i=1}^k$ are eigenvalues of $S^\top S$ and $Y^\top Y$, respectively, that satisfy $0\leq \lambda_i^{S},\lambda_i^Y\leq 1$ and $\lambda_i^{S}+\lambda_i^Y= 1$ for $1\leq i\leq k$. Then
\begin{equation*}
\begin{aligned}
&\|VV^\top - V^*{V^*}^\top\|_F^2
= \operatorname{tr}((VV^\top - V^*{V^*}^\top)^2)=\operatorname{tr}(VV^\top+V^*{V^*}^\top - 2V^*{V^*}^\top VV^\top)\\
&\, = \operatorname{tr}(V^\top V)+\operatorname{tr}({V^*}^\top V^*) - 2\operatorname{tr}(S^\top S) = 2k - 2\operatorname{tr}(S^\top S)   = 2\operatorname{tr}(Y^\top Y)  =2\sum_{i=1}^k\lambda_i^Y.
\end{aligned}
\end{equation*}
Under the assumption $\|VV^\top - V^*{V^*}^\top\|_F \leq 1$, we have
$\sum_{i=1}^k\lambda_i^Y \leq \frac{1}{2}$, which implies
$\lambda_i^S \geq \frac{1}{2}$ for $1\leq  i\leq k$. Therefore
\begin{equation}\label{A2}
A_2\leq -8\gamma\mu\sum_{i=1}^k\lambda_i^S\lambda_i^Y \leq -4\gamma\mu\sum_{i=1}^k\lambda_i^Y
    = -2\gamma\mu\|VV^\top - V^*{V^*}^\top\|_F^2.
\end{equation}
Combining \eqref{A1} and \eqref{A2} with
\eqref{half} to get
\begin{equation*}
\begin{aligned}
\frac{d}{dt}\|VV^\top - V^*{V^*}^\top\|_F^2 & \leq 4\sqrt{k}\,\gamma M\|VV^\top - V^*{V^*}^\top\|_F\|x - x^*\|\\
&\quad-2\gamma\mu\|VV^\top - V^*{V^*}^\top\|_F^2,
\end{aligned}
\end{equation*}
which completes the proof.
\end{proof}

\subsection{Analysis of coupled HiSD}
We analyze the stability of the coupled HiSD (\ref{hisd}) in the following theorem. 
\begin{theorem}\label{i3}
Suppose the Assumption \ref{ass1} holds, $0<\alpha\leq 1$ is selected such that $\hat \mu = \mu-2\alpha L>0$, and 
 \begin{equation}\label{o1}
 \|x(0) - x^* \|  < \min \Big\{\rho, \frac{2 \hat\mu}{M}, \frac{\alpha \mu}{8\sqrt{k} M} \Big\},~ \|V(0)V(0)^\top - V^*{V^*}^\top\|_F<  \frac{\alpha}{2}.
 	 \end{equation}
 Then the following exponential asymptotic stability for $x^*$ and $V^*{V^*}^\top$ holds
 $$ \| x-x^*\|\leq \tilde K(\varepsilon) e^{-\beta(\hat \mu -\varepsilon) t},~\|VV^\top - V^*{V^*}^\top\|_F \leq Q(\varepsilon, \tilde\varepsilon) e^{-\min\{\gamma\mu-\tilde\varepsilon, \, \beta(\hat{\mu} - \varepsilon)\}t},~t\geq 0$$
for some $0<\varepsilon<\hat\mu$, $0 < \tilde\varepsilon < \gamma\mu$ and $Q(\varepsilon, \tilde\varepsilon)>0$.
\end{theorem}
\begin{remark}\label{rem61}
In this theorem, the decay rate of $x-x^*$ is proved to be $\beta(\hat\mu-\varepsilon)=\beta(\mu-2\alpha L-\varepsilon)$, which is independent from the relaxation parameter $\gamma$ in the dynamics of $V$. In other words, the selection of $\gamma$ will not affect the decay rate of $x-x^*$, which is in contrast with the linear stability analysis of the coupled HiSD, where the coupled decay rate of the whole system may be deteriorated by relatively small $\gamma$ and then becomes less than $\beta(\mu-2\alpha L-\varepsilon)$. Since the asymptotic stability of $x^*$ is the main concern, such decoupled estimate  is desirable with sufficient practical significance.
\end{remark}
\begin{proof}
We divide the proof into three steps.

\textbf{Part I: Short-time boundedness.}
By (\ref{o1}) and the continuity of solutions, there exists a $\delta_0>0$ such that for all $t\in [0,\delta_0]$,
$$
x(t)\in B_2:=B\Big(x^*,\min\{\rho,\displaystyle\frac{2 \hat\mu}{M}, \frac{\alpha \mu}{8\sqrt{k} M} \Big\}\Big),~\|VV^\top - V^*{V^*}^\top\|_F \leq \frac{\alpha}{2}.
$$
As $x(t)\in B_2$ for $t\in [0,\delta_0]$, we apply Lemma \ref{C7} to get
$
	\|\hat V \hat V^\top - V^*{V^*}^\top\|_2\leq \frac{\sqrt{2k}M}{\mu}\|x(t)-x^*\|,
$
which implies
\begin{equation}\label{a}
\begin{aligned}
&\|VV^\top-\hat V \hat V^\top\|_2  \leq \|VV^\top-V^* {V^*}^\top\|_2+ \|V^* {V^*}^\top-\hat V \hat V^\top\|_2 \\
 &\qquad \leq \|VV^\top - V^*{V^*}^\top\|_F+\frac{\sqrt{2k}M}{\mu}\|x(t)-x^*\|\leq \frac{\alpha}{2}  + \frac{\alpha}{2}  \leq \alpha.
\end{aligned}
\end{equation}
where we have used $\|x(t)-x^*\|\leq  \frac{\alpha \mu}{2\sqrt{2k} M}$ in the last inequality.

Since $B_2\subset B_0$, we have  $x(t)\in B_0 $ for $0< t\leq \delta_0$. Then, based on \eqref{a}, we following an analogous derivation as the proof of Theorem \ref{i1} to obtain 
\begin{equation*}
 \frac{1}{2}\frac{d}{dt}\|x-x^*\|^2  \leq \beta \Big(-\hat \mu +\frac{M}{2}\|x-x^*\|\Big) \|x-x^*\|^2,
\end{equation*}
and in particular,
\begin{equation*}
\frac{d}{dt}\|x-x^*\| \Big|_{t=0^+}  \leq \beta \Big(-\hat \mu +\frac{M}{2}\|x(0)-x^*\|\Big) \|x(0)-x^*\|<0.
\end{equation*}
Thus, there exists a $0<\delta\leq \delta_0$ such that $\|x-x^*\|<\|x(0)-x^*\|$ for any $0< t\leq \delta$. 

\textbf{Part II: Global boundedness.}
We consider the sets
\begin{equation*}
	U=\{t\geq 0:\|x-x^*\|>\|x(0)-x^*\|\} \text{~and~} H=\{t\geq 0:\|VV^\top - V^*{V^*}^\top\|_F > \alpha\}.
\end{equation*}
We first prove that at least one of $U$ and $H$ is empty. If not, then both $U$ and $H$ are non-empty. Since both sets are bounded below by $0$, their infimums exist and are denoted by $t_0=\inf U$ and $t_1= \inf H$. Since we have shown that $\|x-x^*\|<\|x(0)-x^*\|$ for any $0< t\leq \delta$, it follows that $t_0 \geq\delta>0$. Similarly, since $\|VV^\top - V^*{V^*}^\top\|_F \leq \frac{\alpha}{2}$ for $t\in [0,\delta_0]$, we have $t_1 \geq\delta_0>0$.

Then we find the contradiction case by case:

 \underline{Case 1: $t_1\leq t_0$.} As  $\|VV^\top - V^* {V^*}^\top \|_F\leq  \alpha \leq 1$ and  $x(t) \in B(x^*,\rho)$ for $t\in [0, t_1]$, we apply Lemma \ref{L8}  for $t\in [0,t_1]$ to get 
\begin{equation}\label{ode}
\begin{aligned}
\frac{d}{dt}\|VV^\top - V^*{V^*}^\top\|_F & \leq -\gamma\mu\|VV^\top - V^*{V^*}^\top\|_F
    + 2\sqrt{k}\,\gamma M\|x - x^*\|\\
&\hspace{-0.3in}\leq  -\gamma\mu\|VV^\top - V^*{V^*}^\top\|_F + 2\sqrt{k}\,\gamma M\|x(0) - x^*\|.  
\end{aligned}
\end{equation}
Multiplying both sides of (\ref{ode}) by $e^{\gamma \mu t}$ gives
\begin{equation*}
	\frac{d}{dt}(e^{\gamma \mu t}\cdot\|VV^\top - V^*{V^*}^\top\|_F)\leq 2\gamma \sqrt{k}  M e^{\gamma \mu t}\cdot\|x(0)-x^*\|.
\end{equation*}
Integrating from $0$ to $t\in [0,t_1]$, we obtain
\begin{equation*}
	\begin{aligned}
&\|VV^\top - V^*{V^*}^\top\|_F \\
&\qquad \leq \|V(0)V(0)^\top \!-\! V^*{V^*}^\top\|_F\cdot e^{-\gamma \mu t}  +\frac{2 \sqrt{k} M}{\mu}\|x(0)-x^*\|\cdot(1-e^{-\gamma \mu t})\\
&\qquad \leq \frac{\alpha}{2}+\frac{2 \sqrt{k} M}{\mu}\|x(0)-x^*\|\leq \frac{3\alpha}{4},
	\end{aligned}
\end{equation*}
where we have used $\|x(0)-x^*\|\leq  \frac{\alpha \mu}{8 \sqrt{k} M}$ in the last step. This leads to $\|VV^\top - V^*{V^*}^\top\|_F\Big|_{t= t_1}< \alpha $, which means that $t_1$ can be further enlarged, a contradiction to the definition of $t_{1}$. 

 \underline{Case 2: $t_1>t_0$.} For this case, the $VV^\top $ in the dynamics of $x$ can be viewed as a perturbation of $\hat V\hat V^\top$ such that one can follow the analysis of the perturbed HiSD (i.e. a slight modification of the proof of Theorem \ref{T1}) to reach a contradiction.

Consequently, in both cases we have reached a contradiction such that at least one of $U$ and $H$ is empty. 

We now show that both $U$ and $H$ are empty by excluding other possibilities:
\begin{itemize}
\item If $U$ is empty but $H$ is not empty, we can apply the same argument as Case 1 to show that $t_1$ can be further enlarged, a contradiction to the definition of $t_{1}$.

\item If $H$ is empty but $U$ is not empty,  the $VV^\top $ in the dynamics of $x$ can be again viewed as a perturbation of $\hat V\hat V^\top$ such that one can follow the analysis of the perturbed HiSD (i.e. a slight modification of the proof of Theorem \ref{T1}) to reach a contradiction.
\end{itemize}

In summary, both $U$ and $H$ are empty such that
\begin{equation}\label{zjy}
	\|x-x^*\|\leq \|x(0)-x^*\| \text{~and~}  \|VV^\top - V^*{V^*}^\top\|_F \leq  \alpha,\, \forall t>0.
\end{equation}

\textbf{Part III: Exponential stability.}
With (\ref{zjy}) in hand,  the $VV^\top $ in the dynamics of $x$ can be viewed as a perturbation of $\hat V\hat V^\top$ such that Theorem \ref{i1} gives the exponential asymptotic stability of $x$. Thus, it remains to prove the stability of $V$. Since $\|VV^\top - V^*{V^*}^\top\|_F \leq \alpha \leq 1$ for any $t\geq 0$,  we apply Lemma \ref{L8} and $\|x-x^*\|\leq  \tilde K(\varepsilon) e^{-\beta(\hat \mu -\varepsilon) t}$ (cf. Theorem \ref{i1}) to get
 \begin{equation*}
\begin{aligned}
\frac{d}{dt}\|VV^\top - V^*{V^*}^\top\|_F & \leq -\gamma\mu\|VV^\top - V^*{V^*}^\top\|_F
    + 2\sqrt{k}\,\gamma M\|x - x^*\|\\
&\leq  -\gamma\mu\|VV^\top - V^*{V^*}^\top\|_F + 2\sqrt{k}\,\gamma M\tilde K(\varepsilon) e^{-\beta(\hat \mu -\varepsilon) t}.
\end{aligned}
\end{equation*}
Multiply both sides by  $e^{\gamma \mu t}$ gives
\begin{equation}\label{vv}
	\frac{d}{dt}(e^{\gamma \mu t}\cdot\|VV^\top - V^*{V^*}^\top\|_F)\leq 2\sqrt{k}\, \gamma  M \tilde K(\varepsilon) e^{(\gamma \mu -\beta(\hat \mu -\varepsilon))t}.
\end{equation}
If $\gamma \mu -\beta(\hat \mu -\varepsilon)\neq 0$, integrating this equation from $0$ to $t$ yields
\begin{equation*}
\begin{aligned}
&\|VV^\top - V^*{V^*}^\top\|_2 \leq \|VV^\top - V^*{V^*}^\top\|_F \\
&\qquad \leq \|V(0)V(0)^\top - V^*{V^*}^\top\|_F\cdot e^{-\gamma\mu t} + \frac{2\sqrt{k}\gamma M\tilde{K}(\varepsilon)}{|\gamma\mu - \beta(\hat{\mu} - \varepsilon)|} 
\left|e^{-\beta(\hat{\mu} - \varepsilon)t} - e^{-\gamma\mu t}\right|\\
&\qquad \leq \left(\frac{\alpha}{2}+ \frac{2\sqrt{k}\gamma M\tilde{K}(\varepsilon)}{|\gamma\mu - \beta(\hat{\mu} - \varepsilon)|}\right) e^{-\min\{\gamma\mu, \beta(\hat{\mu} - \varepsilon)\}t}.
\end{aligned}
\end{equation*}
If $\gamma \mu -\beta(\hat \mu -\varepsilon)=0$, then (\ref{vv}) becomes
\begin{equation*}
	\frac{d}{dt}(e^{\gamma \mu t}\cdot\|VV^\top - V^*{V^*}^\top\|_F)\leq 2\sqrt{k}\,  \gamma M \tilde K(\varepsilon),
\end{equation*}
and we integrate this equation from $0$ to $t$ to get
\begin{equation*}
\|VV^\top - V^*{V^*}^\top\|_2 \leq \|VV^\top - V^*{V^*}^\top\|_F  \leq \Big(\frac{\alpha}{2}+ 2\sqrt{k}\,\gamma M\tilde{K}(\varepsilon)t \Big) e^{-\gamma\mu t}.
\end{equation*}
We combine this with $te^{-\tilde\varepsilon t}\leq (e\tilde\varepsilon)^{-1}$ to complete the proof.
\end{proof} 

\section{Stability of momentum-accelerated variants} \label{acce}
In this section, we establish the asymptotic stability of  momentum-accelerated variants of idealized, perturbed, and coupled HiSD. Before proceeding with the analysis, we introduce several preliminary lemmas.

To characterize the sensitivity of eigenvectors under matrix perturbations, we first introduce a variant of the Davis--Kahan theorem.
\begin{lemma}\cite[Corollary 1]{Luoacc} \label{C7}
Suppose the Assumption \ref{ass1} holds. For any $x,y \in B(x^*,\rho)$, consider the orthogonal projections $\sum_{i=1}^k v_{x,i}v_{x,i}^\top $ and $\sum_{i=1}^k v_{y,i}v_{y,i}^\top $ where $\{v_{x,i}\}_{i=1}^k$ and $\{v_{y,i}\}_{i=1}^k$ are two sets of orthonormal eigenvectors corresponding to the first $k$ smallest eigenvalues of $\nabla^2 E(x)$ and $\nabla^2 E(y)$, respectively. Then
\begin{equation*}
	\Big\|\sum_{i=1}^k v_{x,i}v_{x,i}^\top  - \sum_{i=1}^k v_{y,i}v_{y,i}^\top\Big\|_2 \leq \frac{\sqrt{2k}M}{\mu}\|x-y\|_2.
\end{equation*}
\end{lemma}

We then refer the following lemma for future use \cite[Section 29, Theorem V]{Wal}.
 \begin{lemma}\label{lemexp}
 If the real parts of the eigenvalues of a real or complex matrix $A$ are less than $\upsilon$ for some constant $\upsilon$, then the matrix exponential $e^{At}$ could be bounded by $\|e^{At}\|\leq c e^{\upsilon t}$ for some $c\geq 1$
 and for $t\geq 0 $.
 \end{lemma}

\begin{lemma}\label{exp}
For the matrix $A^*$ defined by
\begin{align*}
A^*:=
\begin{bmatrix}
0&I\\
-\alpha_2(I-2V^* {V^*}^\top)\nabla^2 E(x^*)&-\alpha_1 I
\end{bmatrix},
\end{align*}
all its eigenvalues have real parts that are less than or equal to $-\sigma(<0)$ where 
\begin{equation}\label{sigma}
\hspace{-0.1in}-\sigma:=
\left\{
\begin{array}{ll}
\displaystyle\frac{-\alpha_1+\varepsilon}{2},~~&\displaystyle \Delta\leq 0;\\[0.15in]
\displaystyle\frac{-\alpha_1+\sqrt{\alpha_1^2-4\alpha_2\mu}}{2},~~&\displaystyle \Delta> 0,
\end{array}
\right.~\Delta:=\alpha_1^2-4\alpha_2\mu,
 \end{equation}
 for some $0<\varepsilon\ll\alpha_1$. 
\end{lemma}
\begin{proof}
 One could follow a similar spectral analysis as \cite[Theorem 1]{Luoacc} to prove this lemma, and we thus omit the details for simplicity. The difference is that in \cite[Theorem 1]{Luoacc} only the case $\Delta\leq 0$ is considered, and we further consider the case $\Delta>0$ here.
 \end{proof}
 
By the above two lemmas,  there exists a constant $c_0 \geq 1$ such that $\|e^{A^*t}\|\leq c_0 e^{-\sigma t}$ for $t\geq 0$.

\subsection{The idealized case} We consider the asymptotic stability of (\ref{asadk}). Let $y:=[x^\top,m^\top]^\top$,  $y^*:=[(x^*)^\top,0^\top]^\top$ where $x^*$ is the target index-$k$ saddle point of $E$. We first reformulate the second equation of (\ref{asadk}) as 
 \begin{equation}\label{21}
 \begin{aligned}
 \frac{dm}{dt}
 &=-\alpha_1 m-\alpha_2 (I-2\hat V\hat V^\top) (\nabla E(x)-\nabla E(x^*))\\
 &=-\alpha_1 m-\alpha_2(I-2\hat V\hat V^\top)\int_0^1 \nabla^2 E (x^*+s(x-x^*))ds (x-x^*)\\
 &=-\alpha_1 m-\alpha_2\mathcal (I-2 V^* {V^*}^\top)\nabla^2 E (x^*)(x-x^*)+p(x).
 	\end{aligned}
 \end{equation}
Here $V^*:=[v_1^*,\ldots,v_k^*]\in\mathbb R^{N\times k}$ where  $\{v_i^*\}_{i=1}^k$ are orthonormal eigenvectors of $\nabla^2 E(x^*)$ corresponding to its $k$ smallest eigenvalues, and
\begin{equation*}
 p:=  -\alpha_2\Big( (I-2\hat V\hat V^\top) \int_0^1 \nabla^2 E(x^*+s(x-x^*))ds -\mathcal (I-2 V^* {V^*}^\top) \nabla^2 E(x^*)\Big)(x-x^*).
\end{equation*}

 \begin{theorem}\label{aa}
 Suppose the Assumptions \ref{ass1} hold and
\begin{equation}\begin{array}{l}\label{y0}
\displaystyle \|y(0)-y^*\|< \frac{J}{c_0},~
J:=\min\Big\{\rho,\frac{\sigma\theta}{c_0\alpha_2C_1}\Big\},
~C_1:=\frac{2\sqrt{2k}ML}{\mu}+\frac{M}{2},
\end{array}
\end{equation}
for some $0<\theta< 1$. Then the stationary point $y^*$ of the accelerated HiSD (\ref{asadk}) is exponential asymptotic stable and
\begin{align*}
\|y-y^*\|\leq c_0\|y(0)-y^*\|e^{-\sigma(1-\theta) t},~~t\geq 0.
\end{align*}
\end{theorem}
 \begin{proof}
By (\ref{y0}), the continuity of the solutions and $c_0\geq 1$, there exists a $\delta_0>0$ such that $x(t)\in B_1:=B(x^*,J)$ for $t\in [0,\delta_0]$. Then for $t\in [0,\delta_0]$, we apply the Assumptions \ref{ass1} and Lemma \ref{C7} to get 
\begin{align}
\|p\|&=\Big\|-\alpha_2\Big((I-2\hat V\hat V^\top) \int_0^1\big(\nabla^2 E(x^*+s(x-x^*)) \!-\!\nabla^2 E(x^*)\big)ds\label{p}\\
&~~+ ((I-2\hat V\hat V^\top)-(I-2V^* {V^*}^\top))\nabla^2 E(x^*)\Big)(x-x^*)\Big\| \leq \alpha_2C_1\|x-x^*\|^2,\nonumber
\end{align}
where $C_1$ is given in (\ref{y0}). 

Let $p_0=[0, p]^\top$ such that (\ref{21}) and the first equation of (\ref{asadk}) could be reformulated as the matrix form
\begin{equation}\label{ie}
		\frac{d (y-y^*)}{dt}=A^*(y-y^*)+p_0~\Rightarrow~ y-y^*=e^{A^*t}(y(0)-y^*)+\int_0^te^{A^*(t-s)}p_0ds.
\end{equation}
We multiply both sides of (\ref{ie}) by $(y-y^*)^\top$ to get
\begin{equation*}
	\begin{aligned}
		\|y-y^*\|^2= (y-y^*)^\top e^{A^*t}(y(0)-y^*)+\int_0^t (y-y^*)^\top e^{A^*(t-s)}p_0ds.
	\end{aligned}
\end{equation*}
Applying the bound $\|e^{A^*t}\|\leq c_0e^{-\sigma t}$ and the estimate (\ref{p}), and dividing both sides by $\|y-y^*\|$, we obtain
\begin{equation}\label{ls7}
\begin{aligned}
\|y-y^*\| & \leq c_0e^{-\sigma t}\|y(0)-y^*\|+c_0 \int_0^te^{-\sigma(t-s)}\|p_0\|ds \\
&\leq c_0e^{-\sigma t}\|y(0)-y^*\|+c_0\alpha_2C_1 \int_0^te^{-\sigma(t-s)}\|y-y^*\|^2ds.	
\end{aligned}
\end{equation}

Now we consider the set
$W=\{t\geq 0:\|y-y^*\|> J\}$. We intend to prove that $W$ is empty. If not, then $W$ is a non-empty set with a lower bound $0$ and thus has an infimum denoted by $t_0$. We first apply (\ref{y0}) and $c_0\geq 1$ to obtain $\|y(0)-y^*\|\leq c_0 \|y(0)-y^*\|< J$. Thus, there exists a $0<\delta\leq \delta_0$ such that $\|y-y^*\|\leq J$ for any $0\leq t\leq \delta$. For this reason, $t_0$ can not be less than $\delta$ and is thus positive. Furthermore, the estimate (\ref{p}) and thus (\ref{ls7}) remain valid for $t\in [0,t_0]$ due to $x\in B_1$ over  $t\in [0,t_0]$. Then we combine (\ref{ls7}), (\ref{y0}) and the continuity of $y-y^*$ to get a contradiction as follows
  \begin{align*}
J&=\|y(t_0)-y^*\| \leq c_0e^{-\sigma t_0}\|y(0)-y^*\|+c_0\alpha_2C_1 \int_0^{t_0}e^{-\sigma(t_0-s)}\|y-y^*\|^2ds\\
&< e^{-\sigma t_0}J+c_0\alpha_2C_1 J^2\int_0^{t_0}e^{-\sigma(t_0-s)}ds\\
&\leq J\Big(e^{-\sigma t_0}+c_0\alpha_2C_1 \frac{\sigma\theta}{c_0\alpha_2C_1 }\frac{1-e^{-\sigma t_0}}{\sigma}\Big)=J\big(\theta+(1-\theta)e^{-\sigma t_0}\big)<J.
\end{align*}
Consequently, $\|y-y^*\|\leq J$ for any $t\geq 0$, which further implies that (\ref{ls7}) holds for any $t\geq 0$. We apply this to rewrite (\ref{ls7}) in terms of $e^{\sigma t}\|y-y^*\|$ as
$$e^{\sigma t}\|y-y^*\| \leq c_0\|y(0)-y^*\|+c_0\alpha_2C_1 J\int_0^t e^{\sigma s}\|y-y^*\| ds. $$
An application of the Gronwall inequality yields 
$$e^{\sigma t}\|y-y^*\| \leq c_0\|y(0)-y^*\|\exp\big(c_0\alpha_2C_1 Jt\big),$$ 
which, together with the definition of $J$, implies
\begin{align*}
\|y-y^*\| &\leq c_0\|y(0)-y^*\|\exp\big(-(\sigma-c_0\alpha_2C_1 J)t\big)\leq c_0\|y(0)-y^*\|e^{-(\sigma-\sigma\theta)t}.
\end{align*}
Thus we complete the proof.
 \end{proof}
  
\subsection{The perturbed case}\label{pac} We  proceed to prove exponential asymptotic stability of  (\ref{ieasadk}). In general, the proof follows that of Theorem \ref{aa}, while more subtle estimates are required to account for the impact of the perturbed unstable subspace. 
\begin{theorem}\label{i2} Suppose Assumption \ref{ass1} holds, $\alpha$ is selected such that $\sigma\theta -2c_0\alpha \alpha_2 L>0$ for some $0<\theta< 1$, 
\begin{equation}\begin{array}{l}\label{acc}
\displaystyle \|y(0)-y^*\|< \frac{\tilde J}{c_0},
~\tilde J:=\min\Big\{\rho,\frac{\sigma\theta - 2c_0\alpha \alpha_2 L }{c_0\alpha_2C_1}\Big\},~ C_1=\frac{2\sqrt{2k}ML}{\mu}+\frac{M}{2},
\end{array}
\end{equation}
 and $\sigma$ is given by (\ref{sigma}). Then the stationary point $y^*$ of the perturbed accelerated HiSD (\ref{ieasadk}) is exponential asymptotic stable and
\begin{align*}
\|y-y^*\|\leq c_0\|y(0)-y^*\|e^{-\sigma(1-\theta) t},~~t\geq 0.
\end{align*}
\end{theorem}
\begin{remark}
We compare this convergence rate with that for the accelerated HiSD (\ref{asadk}) proved in Theorem \ref{aa}. Both this theorem and Theorem \ref{aa} suggest $e^{-\sigma(1-\theta)t}$ convergence rate. However, the ranges of $\theta$ are different. In Theorem \ref{aa},  $\theta$ can be selected as close as $0$, while in this theorem $\theta$ should be  bigger than $2c_0\alpha \alpha_2 L/\sigma$, which leads to a lower convergence rate. Such difference indicates the impact of the perturbation on the unstable subspace characterized by $\alpha$. 
\end{remark}
\begin{proof}
Similar to (\ref{21}), we perform a reformulation of the second equation of perturbed accelerated HiSD (\ref{ieasadk}) as 
\begin{equation*}
\begin{aligned}
\frac{dm}{dt}
 &=-\alpha_1 m-\alpha_2\mathcal (I-2\tilde V\tilde V^\top)(\nabla E(x)-\nabla E(x^*))\\
 &=-\alpha_1 m-\alpha_2\mathcal (I-2 V^* {V^*}^\top)\nabla^2 E (x^*)(x-x^*)+\tilde p(x),
	\end{aligned}
\end{equation*}
\vspace{-0.1in}
\begin{equation*}
\begin{aligned}
\ds \tilde p(x):= & -\alpha_2\Big((I-2\tilde V \tilde V^\top) \int_0^1 \nabla^2 E(x^*+s(x-x^*))ds\\
&\qquad -(I-2 V^* {V^*}^\top) \nabla^2 E(x^*)\Big)(x-x^*).
	\end{aligned}
\end{equation*}

By (\ref{acc}), the continuity of the solutions and $c_0\geq 1$, there exists a $\delta_0>0$ such that $ x(t)\in B(x^*,\tilde J)$ for $t\in [0,\delta_0]$. Then for $t\in [0,\delta_0]$, we apply the Assumptions \ref{ass1}, Lemma \ref{C7} and $\|\hat V {\hat V}^\top - \tilde{V} \tilde{V}^\top \|_2 \leq  \alpha$ to get 
\begin{equation*}
	\begin{aligned}
&\|\tilde p\| =\Big\|\alpha_2\Big((I-2\tilde V \tilde V^\top) \int_0^1\big(\nabla^2 E(x^*+s(x-x^*)) - \nabla^2 E(x^*)\big)ds \\
& + \big[\big((I\!-\!2\tilde V \tilde V^\top)\!-\!(I\!-\!2\hat V \hat V^\top)\big)+\big((I\!-\!2\hat V \hat V^\top)\!-\!(I\!-\!2V^* {V^*}^\top)\big)\big]\nabla^2 E(x^*)\Big)(x\!-\!x^*)\Big\|\\
&~~\quad \leq 2\alpha \alpha_2 L\|x-x^*\|+\alpha_2C_1\|x-x^*\|^2.
\end{aligned}
\end{equation*}
 Let $\tilde p_0:=[0,~\tilde p]^\top$, and then a similar derivation as (\ref{ie})--(\ref{ls7}) yields 
\begin{align}\label{aey}
\|y-y^*\|&\leq c_0e^{-\sigma t}\|y(0)-y^*\|+c_0 \int_0^te^{-\sigma(t-s)}\|\tilde p_0\|ds \nonumber\\
&\leq c_0e^{-\sigma t}\|y(0)-y^*\|+c_0 \int_0^te^{-\sigma(t-s)}(\alpha_2C_1\|x-x^*\|^2+ 2\alpha \alpha_2 L\|x-x^*\|)ds \nonumber\\
&\leq c_0e^{-\sigma t}\|y(0)-y^*\|+c_0\alpha_2C_1 \int_0^te^{-\sigma(t-s)}\|y-y^*\|^2ds \\
&\qquad + 2c_0\alpha \alpha_2 L \int_0^te^{-\sigma(t-s)}\|y-y^*\|ds.\nonumber
\end{align}

Now we consider the set
$\tilde W=\{t\geq 0: \|y-y^*\|> \tilde J\}$. We intend to prove that $\tilde W$ is empty. If not, then $\tilde W$ is a non-empty set with a lower bound $0$ and has an infimum denoted by $t_0$. We first apply  (\ref{acc}) and $c_0\geq 1$ to obtain $\|y(0)-y^*\|\leq c_0 \|y(0)-y^*\|< \tilde J$. Thus, there exists a $0<\delta\leq \delta_0$ such that $\|y-y^*\|\leq \tilde J$ for any $0\leq t\leq \delta$. For this reason, $t_0$ can not be less than $\delta$ and is thus positive. Furthermore, the estimate (\ref{aey}) remains valid for $t\in [0,t_0]$ due to $\|y-y^*\|\leq \tilde J$ over  $t\in [0,t_0]$. Then we derive a contradiction from (\ref{aey}) as follows
\begin{align*}
\tilde J&=\|y(t_0)-y^*\| \leq c_0e^{-\sigma t_0}\|y(0)-y^*\|+c_0\alpha_2C_1 \int_0^{t_0}e^{-\sigma(t_0-s)}\|y-y^*\|^2ds \\
&\qquad + 2c_0\alpha \alpha_2 L \int_0^{t_0} e^{-\sigma(t_0-s)}\|y-y^*\|ds\\
&< e^{-\sigma t_0}\tilde J+c_0\alpha_2C_1 \tilde J^2\int_0^{t_0}e^{-\sigma(t_0-s)}ds+ 2c_0\alpha \alpha_2 L \tilde J \int_0^{t_0} e^{-\sigma(t_0-s)} ds\\
&\leq \tilde J\Big(e^{-\sigma t_0}+c_0\alpha_2C_1\Big(\frac{\sigma\theta-2c_0\alpha \alpha_2 L }{c_0\alpha_2C_1} \Big)\frac{1-e^{-\sigma t_0}}{\sigma}+ 2c_0\alpha \alpha_2 L \frac{1-e^{-\sigma t_0}}{\sigma} \Big)\\
&=\tilde J\big(\theta+(1-\theta)e^{-\sigma t_0}\big)<\tilde J.
\end{align*}
Consequently, $\|y-y^*\|\leq \tilde J$ for any $t\geq 0$, which further implies that (\ref{aey}) holds for any $t\geq 0$. Thus we  rewrite (\ref{aey}) in terms of $e^{\sigma t}\|y-y^*\|$ as
$$e^{\sigma t}\|y-y^*\| \leq c_0\|y(0)-y^*\|+(c_0\alpha_2C_1\tilde J+ 2c_0\alpha \alpha_2 L)\int_0^t e^{\sigma s}\|y-y^*\|ds. $$
An application of the Gronwall inequality yields
$$e^{\sigma t}\|y-y^*\| \leq c_0\|y(0)-y^*\|\exp\big((c_0\alpha_2C_1\tilde J+ 2c_0\alpha \alpha_2 L)t\big), $$
which, together with the definition of $\tilde J$, implies
\begin{align*}
\|y-y^*\| \leq c_0\|y(0)\!-\!y^*\|\exp\big(((c_0\alpha_2C_1\tilde J+ 2c_0\alpha \alpha_2 L)\!-\!\sigma)t\big)\leq c_0\|y(0)\!-\!y^*\|e^{-(\sigma-\sigma\theta)t}.
\end{align*}
Thus we complete the proof.
\end{proof}
\begin{remark}
When \( \alpha = 0 \), which corresponds to the cases with exact eigenvectors, Theorem \ref{i1} reduces to Theorem \ref{T1} and Theorem \ref{i2} reduces to Theorem \ref{aa}, which indicates the compatibility of the analysis. 
\end{remark} 

\subsection{The coupled case}
Following the same procedure as in \eqref{21}, we reformulate the second equation of  \eqref{ieah} as 
\begin{equation*}
\begin{aligned}
\frac{dm}{dt} &=-\alpha_1 m-\alpha_2\mathcal (I-2 V^* {V^*}^\top)\nabla^2 E (x^*)(x-x^*)+\bar p(x),
\end{aligned}
\end{equation*}
\vspace{-0.2in}
\begin{equation*}
\begin{aligned}
\ds \bar p(x)&:= \alpha_2((I-2 V^* {V^*}^\top) \nabla^2 E(x^*)\!-\!(I-2 V V^\top) \int_0^1 \nabla^2 E(x^*\!+\!s(x-x^*))ds )(x-x^*)\\
&=\alpha_2\big(\big[(I-2V^* {V^*}^\top)- (I-2V{V}^\top)\big]\nabla^2 E(x^*)\\
&\qquad\qquad  -(I-2 V V^\top) \int_0^1\big(\nabla^2 E(x^*+s(x-x^*)) - \nabla^2 E(x^*)\big)ds\big)(x-x^*),
	\end{aligned}
\end{equation*}
where the second right-hand side term of $\bar p$ can be bounded by the local Lipschitz continuity assumption like (\ref{p}). Compared with the analysis of momentum-accelerated perturbed HiSD in Section \ref{pac}, estimating the first right-hand side term of $\bar p$ in the coupled case additionally requires an a priori bound on $\|VV^\top-V^*{V^*}^\top\|_F$, which can be given by Lemma \ref{L8}. The stability analysis of \eqref{ieah} then follows the same strategy as the proof of Theorem \ref{i3}, with the estimate for $\|x-x^*\|$ obtained analogously to that in Theorem \ref{i2}. We therefore omit the proof and state the results.
\begin{theorem}\label{i4}
Suppose the Assumption \ref{ass1} holds and $0<\alpha\leq 1$ is selected such that $\sigma\theta -2c_0\alpha \alpha_2 L>0$ for some $0<\theta< 1$. Let $\bar J:=\min\Big\{\rho,\frac{2\sigma\theta - 4c_0\alpha \alpha_2 L }{c_0\alpha_2 M}\Big\}$. If 
\begin{equation*}
\displaystyle \|y(0)-y^*\|<\min \Big\{\rho,\, \frac{\bar J}{c_0},\,\frac{\alpha \mu}{8\sqrt{k} M}\Big\} , ~\|V(0)V(0)^\top - V^*{V^*}^\top\|_F < \frac{\alpha}{2},
\end{equation*}
then the following exponential asymptotic stability for $x^*$ and $V^*{V^*}^\top$  of (\ref{ieah}) holds
 $$\|y-y^*\|\leq c_0\|y(0)-y^*\|e^{-\sigma(1-\theta) t}, ~\|VV^\top - V^*{V^*}^\top\|_F \leq \bar Q(\tilde\varepsilon)e^{-\min\{\gamma\mu-\tilde\varepsilon, \, \sigma(1-\theta)\}t},$$
for $t\geq 0$, where $\sigma$ is given by (\ref{sigma}), $0 < \tilde\varepsilon < \gamma\mu$, and $\bar Q(\tilde\varepsilon)>0$.
\end{theorem}

\section{Numerical experiments}\label{num}
We present numerical experiments to substantiate theoretical results.
\subsection{Selection of relaxation parameters}\label{sec61}
We intend to determine relaxation parameters such that the accelerated HiSD converges faster than the idealized HiSD. Equivalently, according to Theorems \ref{T1} and \ref{aa},  we will determine $(\alpha_1,\alpha_2)$ in (\ref{asadk}) such that $\beta(\mu-\varepsilon)<\sigma(1-\theta)$ for a given $\beta$ in (\ref{sadk}). 

\begin{figure}[H]
  	\setlength{\abovecaptionskip}{0pt}
  	\centering
  	\,\includegraphics[width=5.1in,height=2.5in]{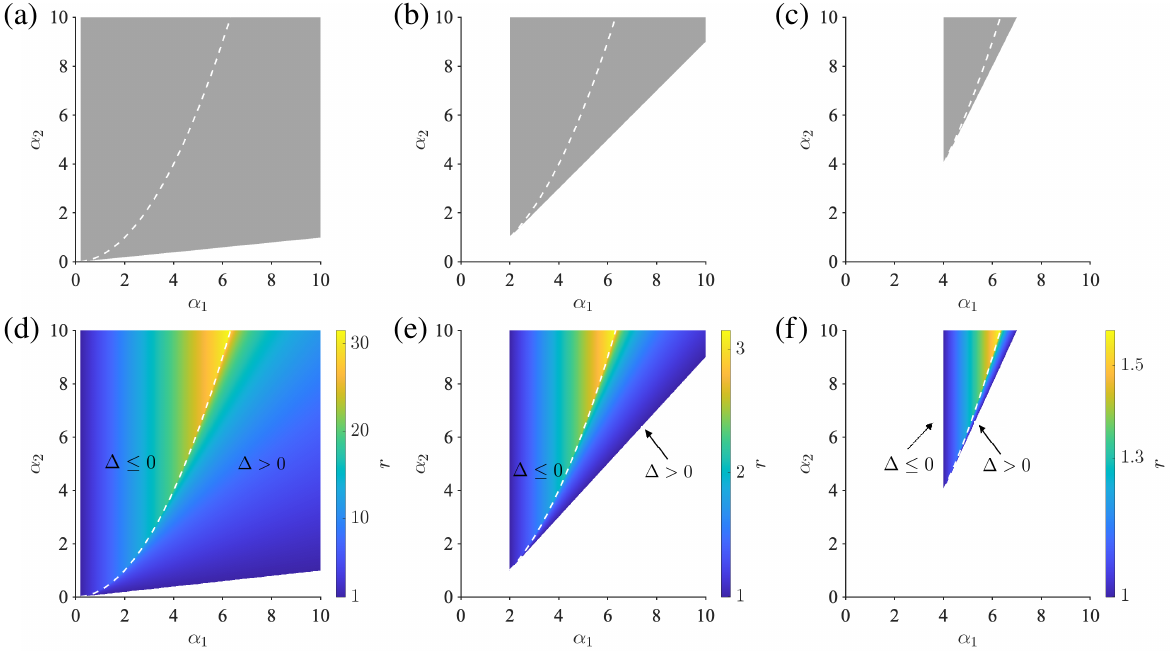}
 \caption{Feasible regions for $(\alpha_1, \alpha_2)$ satisfying 
$\beta(\mu-\varepsilon) < \sigma(1-\theta)$ with $\beta = 0.1$ (left), 
$1$ (middle), and $2$ (right). 
(a)--(c): regions identified numerically; 
(d)--(f): regions determined by \eqref{suf}. 
The dashed curves indicate $\Delta:=\alpha_1^2-4\alpha_2 \mu= 0$. } 
  	\label{region1}
\end{figure}

Since  $\alpha_1$ and $\alpha_2$ are coupled in $\sigma$, it is difficult to derive a necessary and sufficient condition under which the condition $\beta(\mu-\varepsilon)<\sigma(1-\theta)$ holds. Thus, we first perform a numerical investigation. Specifically, we fix $\mu=1$ and $\theta=\varepsilon=0.01$, and consider $\beta=0.1, 1, 2$. For $(\alpha_1,\alpha_2)\in[0,10]^2$, a grid search is performed on a $1000\times1000$ uniform mesh, and the points satisfying $\beta(\mu-\varepsilon)<\sigma(1-\theta)$ are marked. The resulting feasible regions are shown as the shaded areas in Figure \ref{region1} (a)--(c), which indicate that a smaller relaxation parameter $\beta$ leads to a larger feasible region for $(\alpha_1,\alpha_2)$, thereby offering greater flexibility in relaxation parameter selection. As $\beta$ grows, the feasible region shrinks such that we need to be more careful on selection of relaxation parameters $(\alpha_1,\alpha_2)$ to ensure acceleration. 
\vspace{-0.1in}
\vspace{-0.3in}
To explore more insights, we find a sufficient condition of the condition $\beta(\mu-\varepsilon)<\sigma(1-\theta)$:
\begin{equation}\label{suf}
\left\{
\begin{array}{ll}
\displaystyle \alpha_1>2\Theta+\varepsilon, ~~&\displaystyle \Delta\leq 0;\\[0.05in]
\displaystyle \alpha_1>2\Theta,~ \alpha_2 > \frac{\alpha_1 \Theta-\Theta^2}{\mu}, ~~&\displaystyle \Delta> 0,
\end{array}
\right.~~~\Theta:=\frac{\beta(\mu-\varepsilon)}{1-\theta}.
 \end{equation}
 The region of $(\alpha_1,\alpha_2)$ characterized by this sufficient condition is plotted in Figure \ref{region1} (d)--(f), which agrees well with the numerically identified region, thereby providing a practical and explicit criterion for the selection of $(\alpha_1,\alpha_2)$. 

To further quantify the acceleration effect, we evaluate the ratio $r:=\frac{\sigma (1-\theta)}{\beta(\mu-\varepsilon)}$ as a measure of the improvement on the convergence rate in Figure \ref{region1} (d)--(f), which suggests that as $\beta$ increases, the feasible region that satisfies $ \Delta > 0$ shrinks rapidly compared with the one characterized by $ \Delta \leq 0$. Furthermore, in the feasible  region  characterized by $ \Delta \leq 0$, the values of $r$ are larger than those in the feasible  region  characterized by $ \Delta >0$ when the values of $(\alpha_1,\alpha_2)$ in these two feasible  regions are comparable. All these observations suggest that it is more convenient to achieve a substantial acceleration effect in the feasible  region  characterized by $ \Delta \leq 0$.

\begin{figure}[H]
  	\setlength{\abovecaptionskip}{0pt}
  	\centering
 \includegraphics[width=2.5in,height=1.6in]{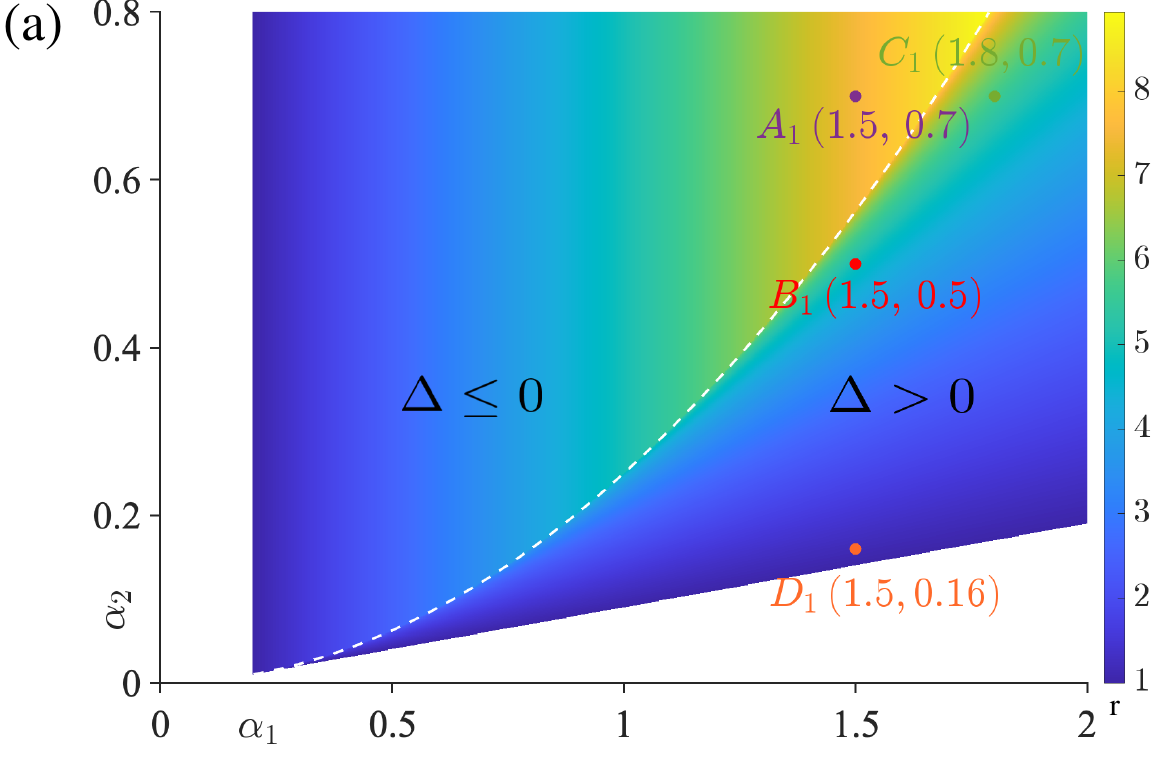}
 \includegraphics[width=2.5in,height=1.6in]{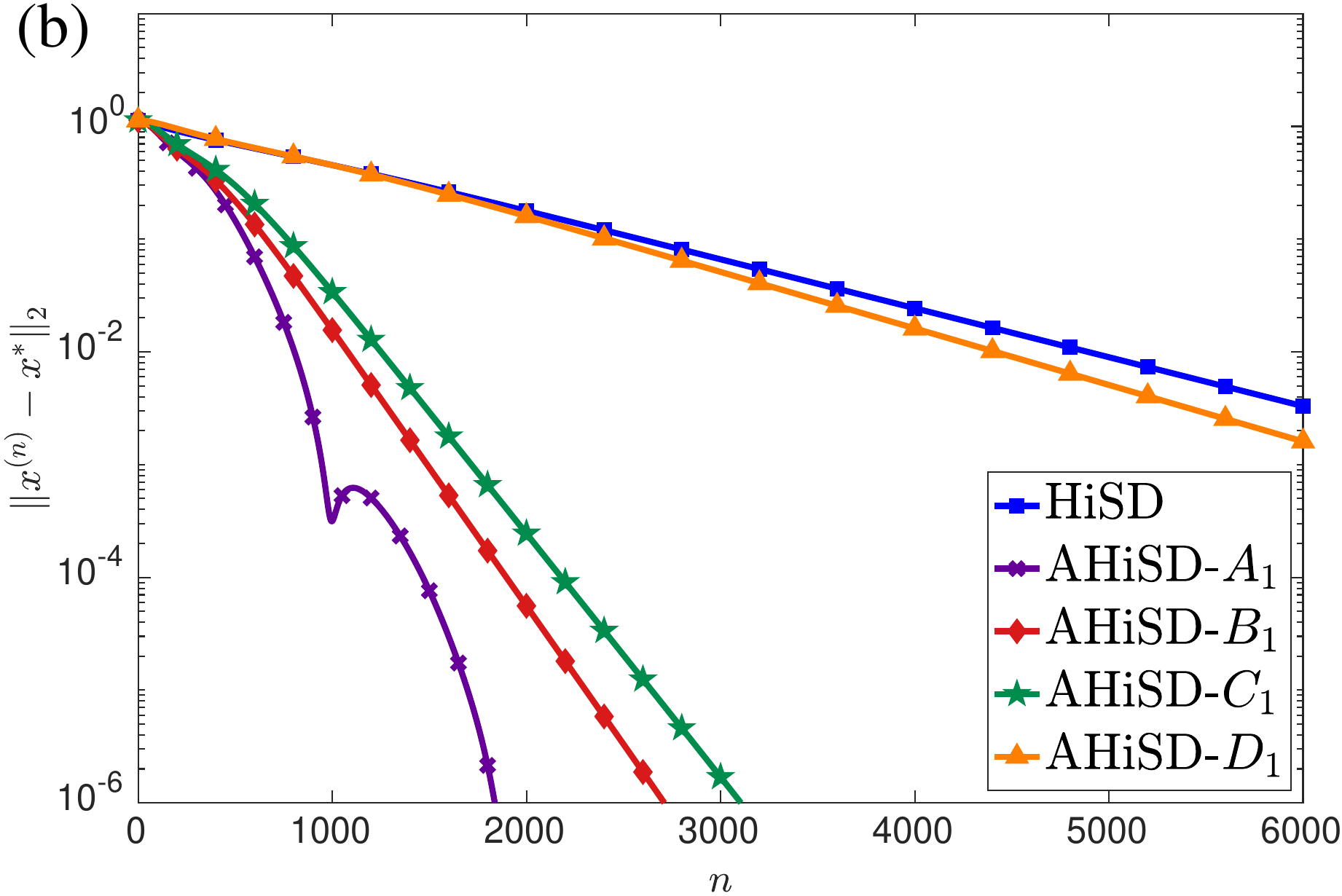}\\[0.05in]
 ~\includegraphics[width=2.5in,height=1.6in]{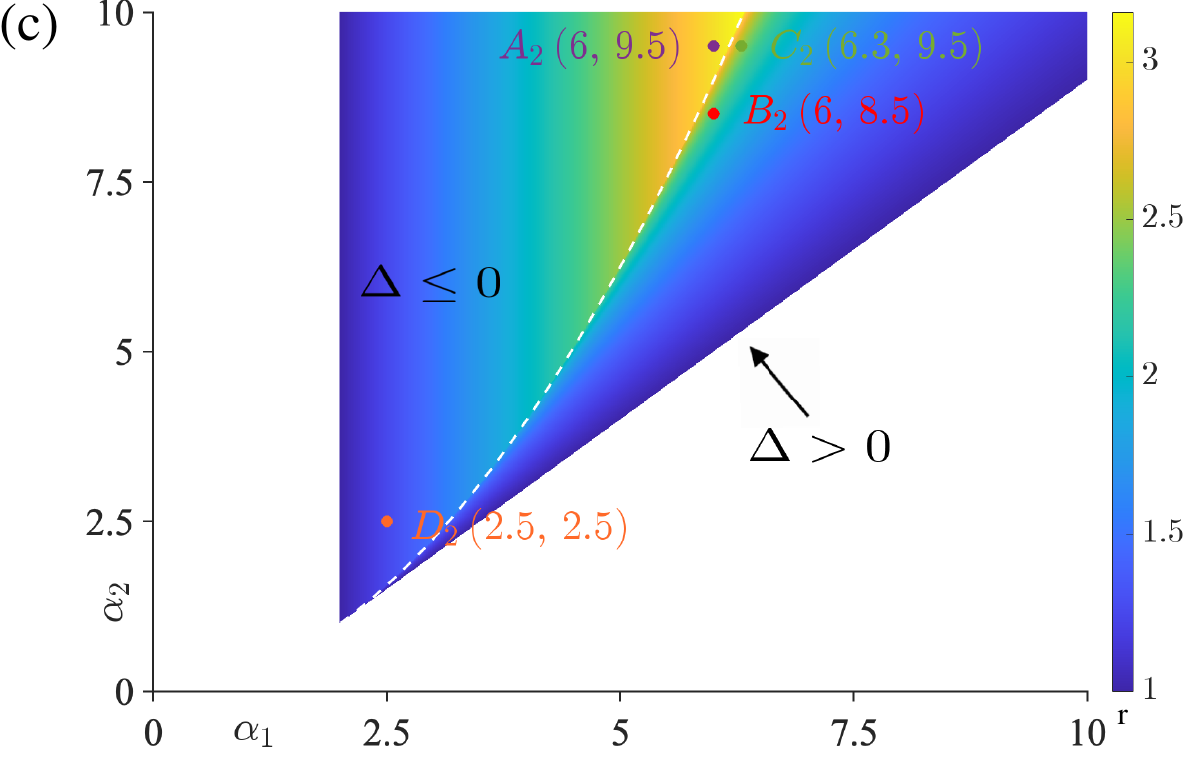}
 \includegraphics[width=2.5in,height=1.6in]{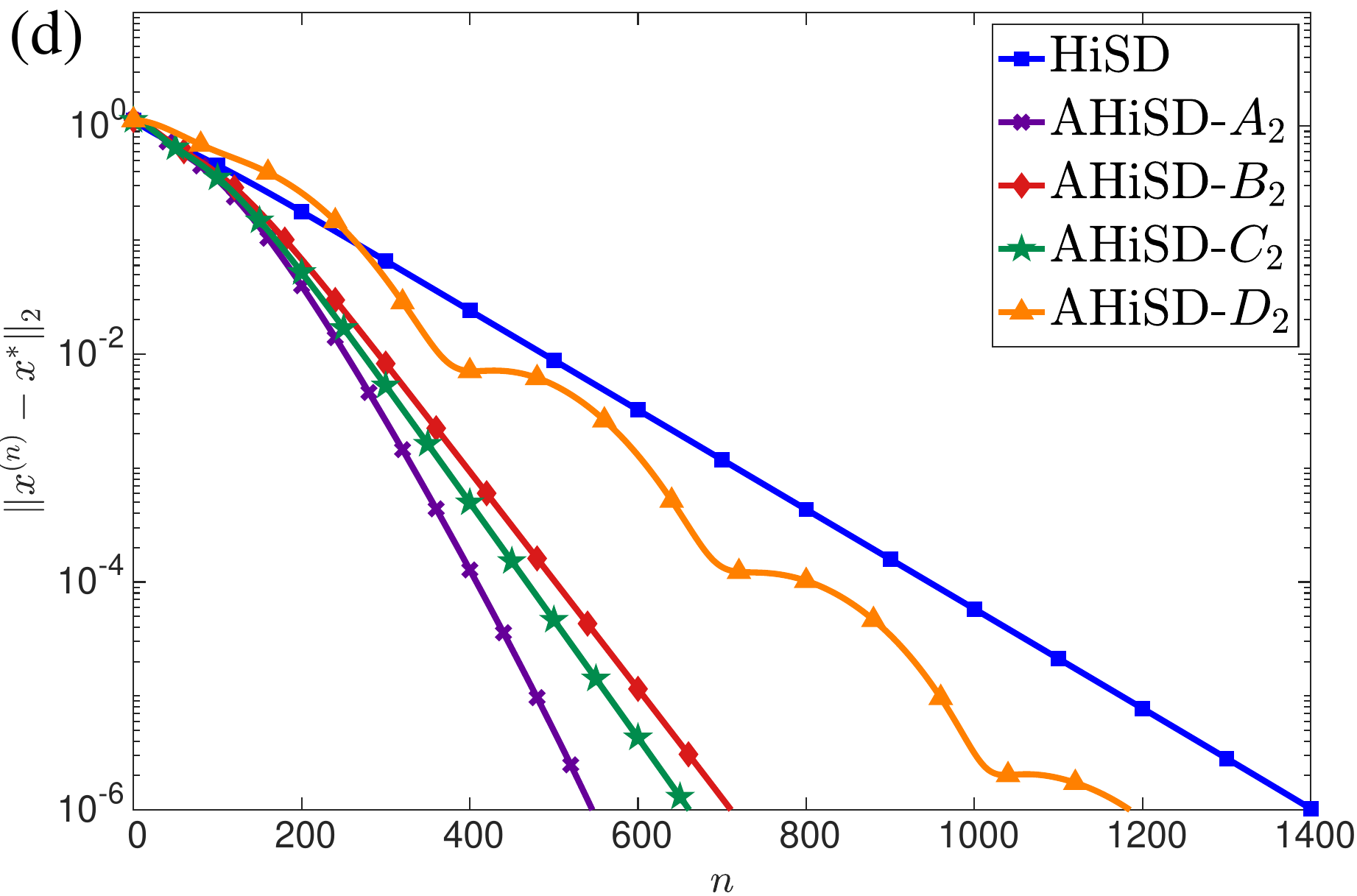}
 \caption{Feasible regions with representative choices of $(\alpha_1, \alpha_2)$ and 
convergence of $\|x^{(n)}- x^*\|_2$ under (a)--(b): $\beta=0.1$ and (c)--(d): $\beta=1$.} 
  	\label{trace}
\end{figure}

To numerically validate the sufficient condition (\ref{suf}), we follow \cite{ZhaSINUM} to explicitly discretize the idealized and the accelerated HiSD as follows
\begin{align}
&\qquad\qquad\frac{x^{(n)}-x^{(n-1)}}{\tau}  = -\beta\bigg(I-2 \hat V^{(n-1)} \hat V^{(n-1)^\top}  \bigg)\nabla E(x^{(n-1)}),\label{back}\\
&\frac{x^{(n+1)}-2x^{(n)}+x^{(n-1)}}{\tau^2}  = -\alpha_1\frac{x^{(n)}-x^{(n-1)}}{\tau} - \alpha_2\bigg(I-2 \hat V^{(n)} \hat V^{(n)^\top}  \bigg)\nabla E(x^{(n)}),\nonumber
\end{align}
where $n\geq 0$, $\tau>0$ is the time step size,
$x^{(-1)}=x^{(0)}$ and the eigenvectors in $\hat V^{(n)}$ are solved by the  simultaneous Rayleigh-quotient iterative minimization method \cite{SRQ}. The $\tau$ is selected to be sufficient small such that the schemes in (\ref{back}) provide good approximations to their continuous counterparts and we remain to investigate the impacts of relaxation parameters $\beta$, $\alpha_1$ and $\alpha_2$. We consider the two-dimensional system $E(x, y)=\frac{1}{4}(x^2-1)^2+ y^2$ proposed in \cite{EZho},
of which $(0,0)$ is an index--$1$ saddle point. The initial point is selected as $ (0.8, 0.8)$, and we fix parameters $\theta = \varepsilon = 10^{-4}$ and the time step size $\tau = 10^{-2}$. The eigenvalues of $\nabla^2 E$ at $(0,0)$ are $-1$ and $2$. We test two cases: $\beta = 0.1$ and $\beta = 1$. For each $\beta$, four representative pairs of $(\alpha_1,\alpha_2)$ are selected, see Figure \ref{trace} for details.

\vspace{-0.3in}

As illustrated in Figure~\ref{trace}, the accelerated HiSD (denoted by AHiSD in the figure) converges faster than the idealized HiSD, and the acceleration effect diminishes as $(\alpha_1, \alpha_2)$ approaches the boundary of the feasible region. For comparable values of $(\alpha_1, \alpha_2)$, parameters in the feasible region characterized by $\Delta \leq 0$ yield better acceleration effect compared to those in the feasible region characterized by $\Delta > 0$.  These observations suggest that it is better to choose $(\alpha_1, \alpha_2)$ within the feasible region characterized by $\Delta \leq 0$ and away from its boundary to achieve satisfactory acceleration effect.

\subsection{Effects of relaxation parameters}\label{sec62} Motivated by  Theorem \ref{i3} and Remark {\ref{rem61}, we investigate the effects of the relaxation parameters $\beta$ and $\gamma$ on the convergence rate.
We consider the three-dimensional energy functional $E(x)=(x_1^2-1)^2 + (x_2^2-1)^2+(x_3^2-1)^2$, which possesses an index-$1$ saddle point $x^*_{1}=(0,\,1,\,1)^\top$ and an index-$2$ saddle point $x^*_{2} =(0,\,0,\,1)^\top$. We apply the explicit scheme in \cite{ZhaSINUM} to discretize (\ref{hisd}), and the initial values are chosen as perturbations of the target saddle point $x^*$ and the corresponding unstable space $V^*$, that is, we take the initial values of $x$ and $V$ as $ x^* + (0.5,\,0.5,\,-0.5)^\top$ and the orthogonal matrix generated by orthonormalizing the columns of $V^* + 0.05\,\boldsymbol{1}_{3\times k}$ via Gram-Schmidt procedure, respectively, where $\boldsymbol{1}_{3\times k}$ is the $3\times k$ all-ones matrix. Throughout the experiments, we set the time-step $\tau=10^{-3}$, and the scheme is terminated when both $\|x-x^*\|\leq 10^{-4}$ and $ \|VV^\top - V^*{V^*}^\top\|_F\leq 10^{-4}$ are satisfied.
\begin{figure}[H]
  	\setlength{\abovecaptionskip}{0pt}
  	\centering
  \includegraphics[width=2.5in,height=1.7in]{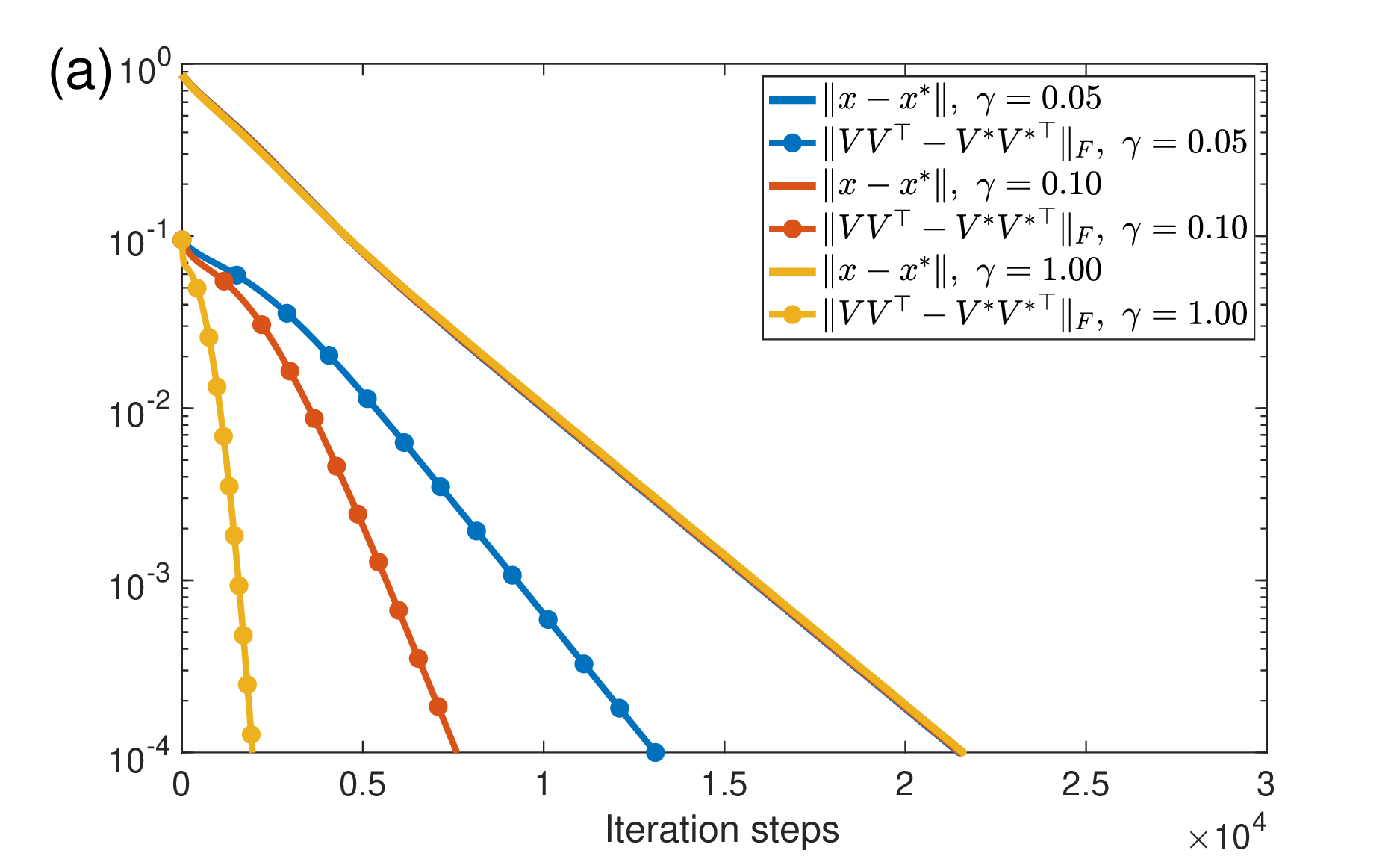}\!\! \includegraphics[width=2.5in,height=1.7in]{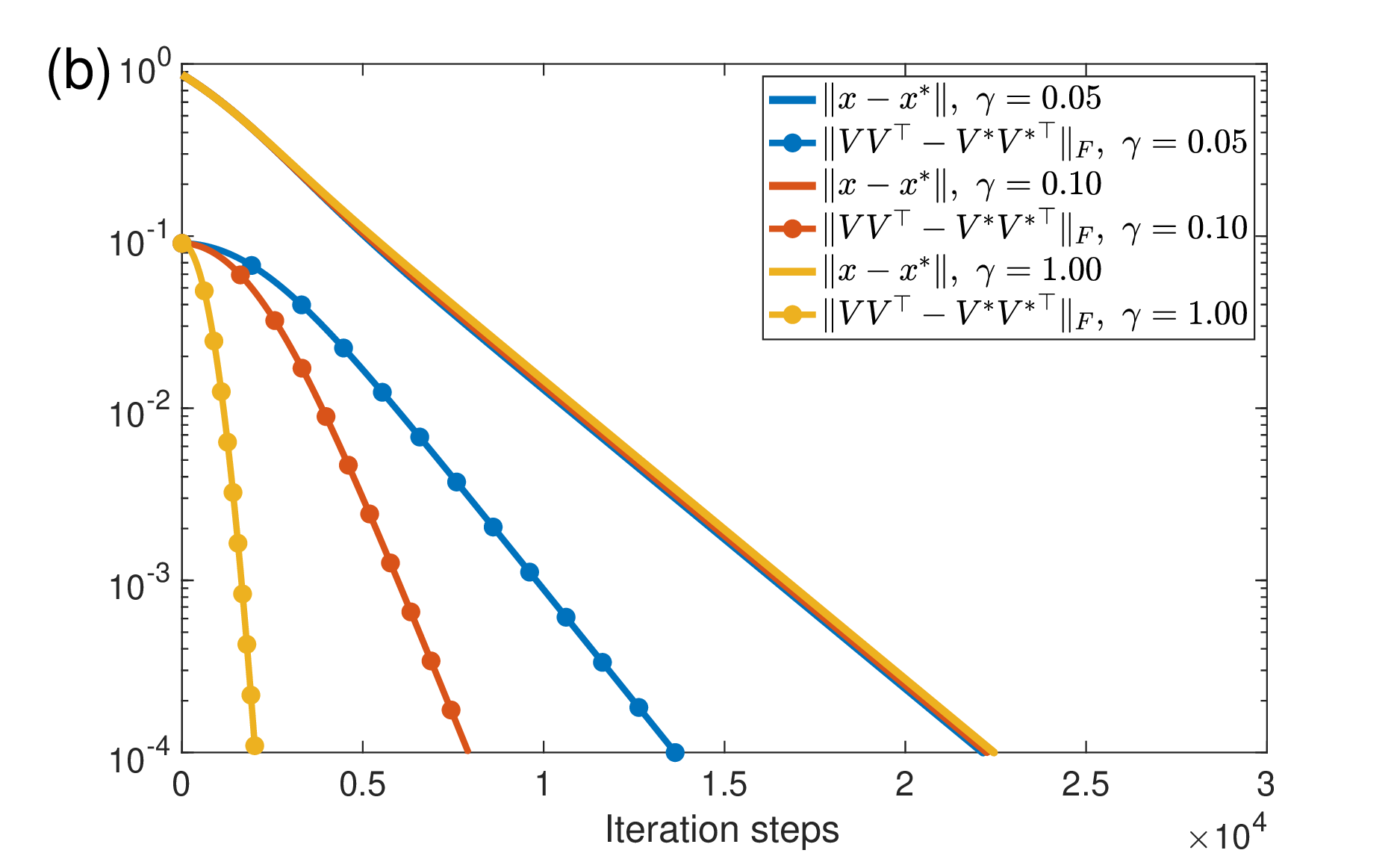}\\
   \includegraphics[width=2.5in,height=1.7in]{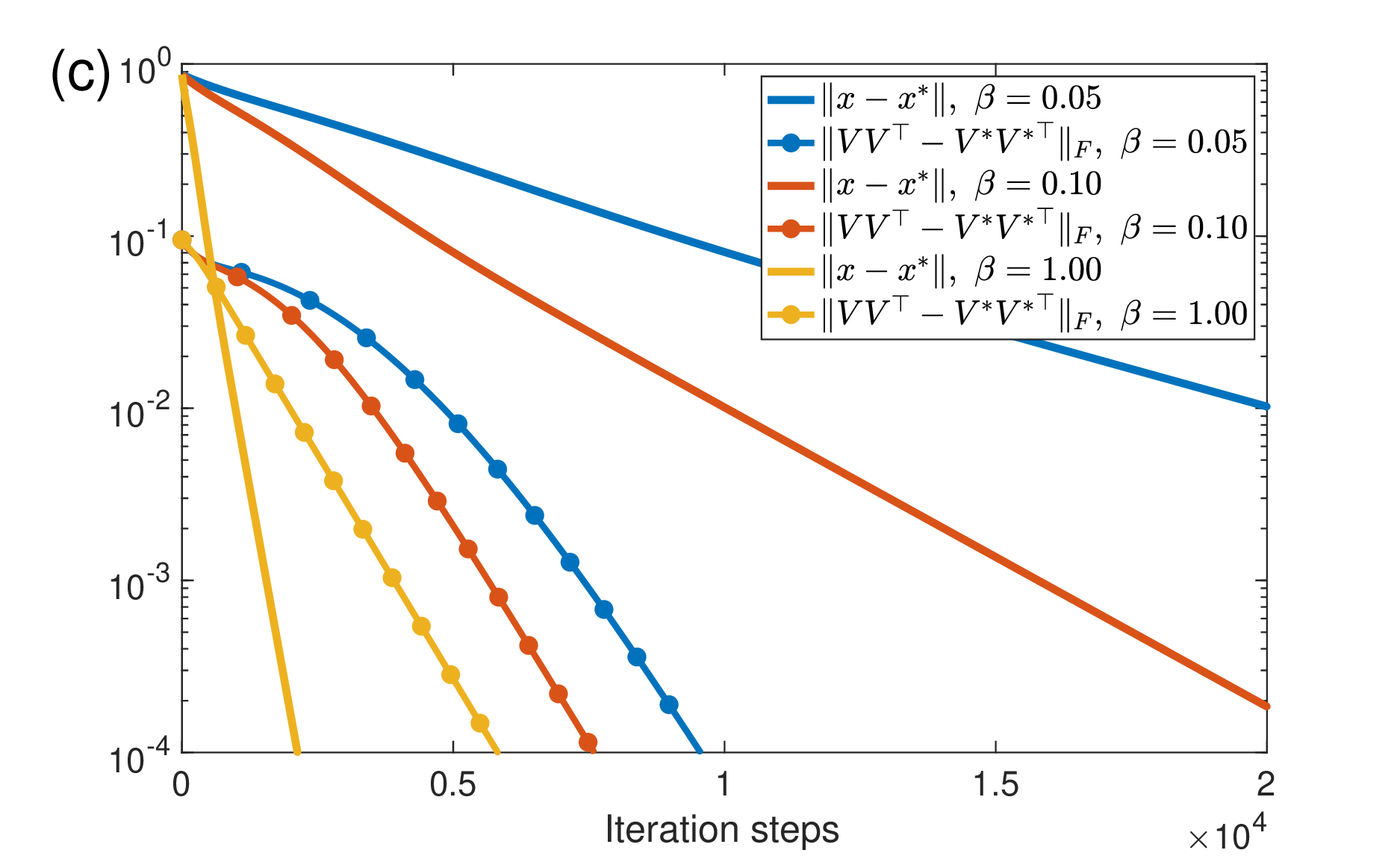}\!\! \includegraphics[width=2.5in,height=1.7in]{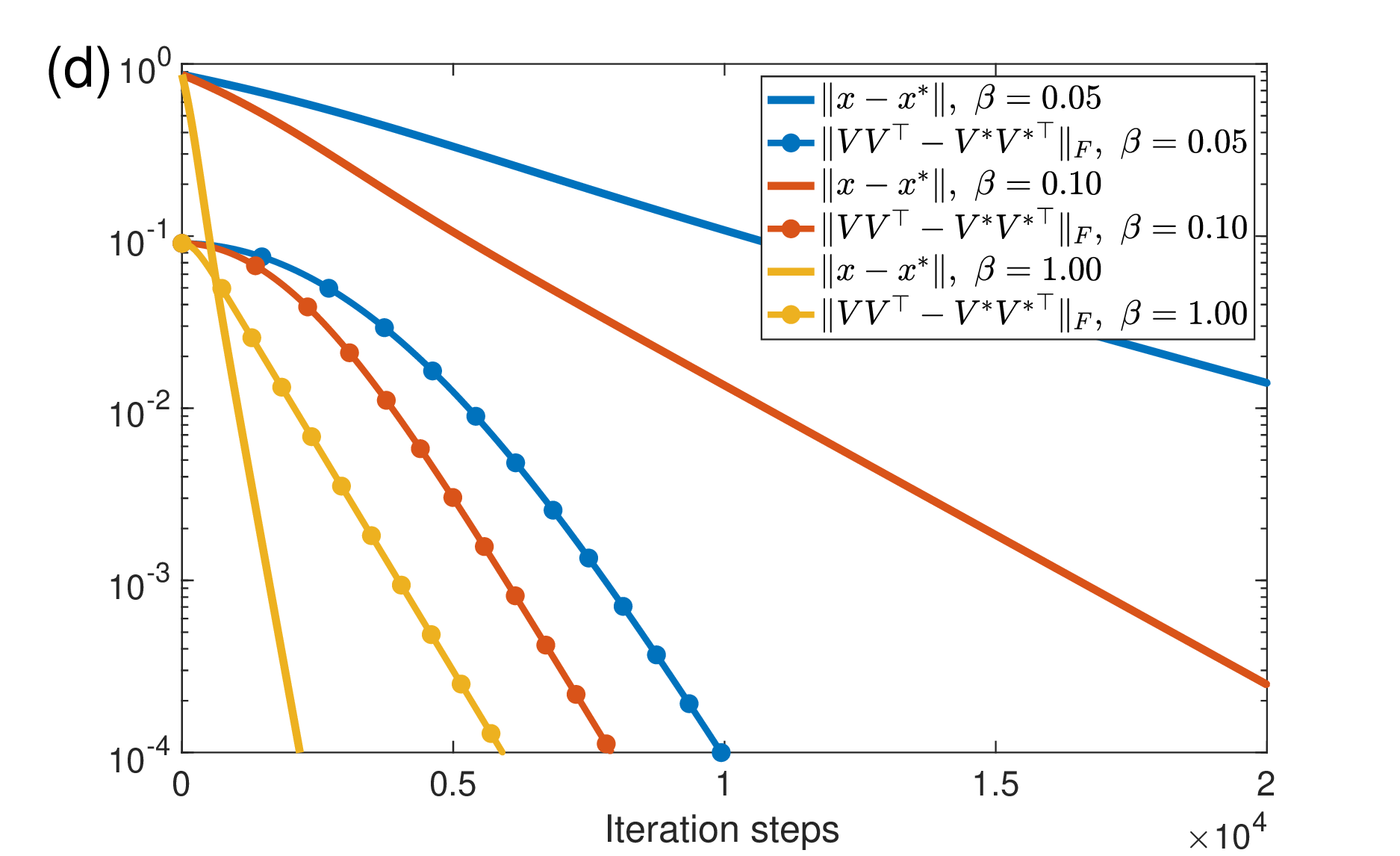}
 \caption{Evolution of  $\|x-x^*\|$ and $\|VV^\top-V^*{V^*}^\top\|_{F}$. The $\beta=0.1$ in   (a)--(b) and the $\gamma=0.1$ in (c)--(d). The saddle point is $x^*_1$ in (a) and (c), and is $x^*_2$ in (b) and (d).} 
  	\label{gamma}
\end{figure}
\vspace{-0.3in}

In (a)--(b) of Figure \ref{gamma}, we fix $\beta=0.1$ to test the convergence behaviors of the coupled HiSD at $x_1^*$ and $x_2^*$. We observe that the value of $\gamma$ has a salient impact on the convergence of the dynamics of $V$, while the convergence of $x$ remains nearly unchanged. These are in good agreement with the theoretical result in Theorem~\ref{i3} that the local convergence rate of the state variable does not depend on that of the unstable subspace.

In (c)--(d) of Figure \ref{gamma}, we fix $\gamma=0.1$ to test the convergence behaviors of the coupled HiSD at $x_1^*$ and $x_2^*$. We observe that when $\beta=0.05$, the dynamics of both $x$ and  $V$ exhibit a slow asymptotic convergence behavior with almost the same convergence rate. These in good agreement with the theoretical results in Theorem~\ref{i3} that the local convergence rate of the dynamics of the unstable subspace can be limited by that of the state variable and in this case both variables have the same predicted asymptotic decay rate. As $\beta$ increases, the convergence of both variables is accelerated.

\section{Concluding remarks}\label{section5}
This work develops asymptotic stability of various HiSD methods, which serves to complement existing stability studies of HiSD.  There are possible directions for future research. For instance, in most existing works the HiSD is proposed under the finite-dimensional setting, while in recent works (cf. \cite{ZhaZhe}) an infinite-dimensional HiSD is proposed for locating saddle points of some semilinear elliptic problems. Thus it is natural to consider the asymptotic stability of such infinite-dimensional HiSD. It is worth mentioning that for the case of partial differential equations,  the linear stability could not imply or exclude the original instability, see e.g. \cite{Kev}. Thus, it is necessary to consider the original asymptotic stability of the infinite-dimensional HiSD for practical applications. We will improve the developed analysis methods to investigate this interesting topic in the near future.


\begin{thebibliography}{10}
\bibitem{Bao} W. Bao and Y. Cai, Mathematical theory and numerical methods for Bose-Einstein condensation. {\it  Kinet. Relat. Models}, 6 (2013), 1--135.

\bibitem{Bof} J. Bofill, W. Quapp,  E.Bernuz, Some remarks on the model of the extended gentlest ascent dynamics. {\it J. Math. Chem.} 53 (2015), 41--57.

\bibitem{Che} H. Chen, G. Dong, J. Iglesias, W. Liu, Z. Xie, Second--order flows for approaching stationary points of a class of nonconvex energies via convex--splitting schemes.
{\it SIAM J. Sci. Comput.} 47 (2025), A1604--A1627.

\bibitem{Chenliu} Z. Chen, W. Liu, Z. Xie, and W. Yi, Nehari manifold optimization and its application for finding unstable solutions of semilinear elliptic PDEs. {\it  SIAM J. Sci. Comput.}  47 (2025), A2098--A2126.

\bibitem{Chu} M. Chu and M. Lin, Generalized gentlest ascent dynamics methods for high-index saddle points. {\it SIAM J. Numer. Anal.} 63 (2025), 2343--2370.

\bibitem{Doye} J. Doye and D. Wales, Saddle points and dynamics of Lennard-Jones clusters, solids, and supercooled liquids. {\it J. Chem. Phys.} 116 (2002), 3777--3788.



\bibitem{du2026convergence}
 Q.~Du and B.~Shi,  Convergence analysis of the discrete constrained
  saddle dynamics and their momentum variants, arXiv preprint
  arXiv:2601.22341,  to appear in {\it Math. Comp.}.



\bibitem{EZho} W. E and X. Zhou, The gentlest ascent dynamics. {\it Nonlinearity} 24 (2011), 1831--1842.


\bibitem{Farr} P. E. Farrell, \'{A}. Birkisson, and S. W. Funke, Deflation techniques for finding distinct solutions of nonlinear partial differential equations. {\it SIAM J. Sci. Comput.} 37 (2015), A2026--A2045.

\bibitem{Farr2} P. E. Farrell, Computing multiple solutions of systems of nonlinear equations
with deflation. In: Proceedings of the International Congress of Mathematicians (2026), Vol. 7. Philadelphia, USA: Society for Industrial and Applied Mathematics and International Mathematical Union.


\bibitem{Gao} W. Gao, J. Leng, and X. Zhou, An iterative minimization formulation for saddle point search. {\it SIAM J. Numer. Anal.} 53 (2015), 1786--1805.

\bibitem{Gou} N. Gould, C. Ortner and D. Packwood, A dimer-type saddle search algorithm with preconditioning and linesearch. {\it Math. Comp.} 85 (2016), 2939--2966.






\bibitem{Kev} P. Kevrekidis, D. Pelinovsky, A. Saxena, When linear stability does not exclude nonlinear instability. {\it Phys. Rev. Lett.} 114 (2015), 214101.







 \bibitem{Li2001} Y. Li and J. Zhou, A minimax method for finding multiple critical points and its applications to semilinear PDEs. {\it SIAM J. Sci. Comput.} 23 (2001), 840--865.
 
\bibitem{SRQ} D. Longsine and S. McCormick, Simultaneous Rayleigh-quotient minimization methods for $Ax=\lambda Bx$. {\it Linear Algebra Appl.}  34 (1980), 195--234.


\bibitem{Luo} Y. Luo, X. Zheng, X. Cheng, L. Zhang, Convergence analysis for discrete high-index saddle dynamics. {\it SIAM J. Numer. Anal.} 60 (2022), 2731--2750.

\bibitem{Luoacc} Y. Luo, L. Zhang, X. Zheng, Accelerated high-index saddle dynamics method for searching high-index saddle points. {\it J. Sci. Comput.} 102 (2025), 31.

\bibitem{Chaves} D. Oyarz\'un, M. Chaves, M. Hoff-Hoffmeyer-Zlotnik, Multistability and oscillations in genetic control of metabolism. {\it J. Theoret. Biol.}  295 (2012), 139--153.


\bibitem{Milnor} J. W. Milnor, {\em Morse Theory}, Princeton University Press, 1963.






\bibitem{Qua} W. Quapp and J. Bofill, Locating saddle points of any index on potential energy surfaces by the generalized gentlest ascent dynamics. {\it Theor. Chem. Acc.} 133 (2014), 1510.

\bibitem{SanWan} E. Sander and T. Wanner, Validated saddle-node bifurcations and applications to lattice dynamical systems. {\it SIAM J. Appl. Dyn. Syst.} 15 (2016), 1690--1733.



\bibitem{Wal} W. Walter, {\it Ordinary differential equations}. Graduate Texts in Mathematics 182, Springer Science and Media, 1998.

\bibitem{wang2021modeling} W. Wang, L. Zhang, P. Zhang, Modelling and computation of liquid crystals. {\it Acta Numer.} 30 (2021), 765--851.

\bibitem{Wu} S. Wu, Y. Xie, D. Pan, L. Zhang, Y. Jin, Hierarchical and ultrametric barriers in the energy landscape of jammed granular matter. {\it Phys. Rev. Lett.,} 136 (2026), p. 198202.








\bibitem{YinPNAS} J. Yin, K. Jiang, A. C. Shi, P. Zhang, L. Zhang, {Transition pathways connecting crystals and quasicrystals}, {\it Proc. Natl. Acad. Sci. U.S.A.} 118 (2021), e2106230118.

%

\bibitem{YinSISC} J. Yin, L. Zhang, P. Zhang, High-index optimization-based shrinking dimer method for finding high-index saddle points. {\it SIAM J. Sci. Comput.} 41 (2019), A3576--A3595.



\bibitem{ZhaDu} J. Zhang, Q. Du, Shrinking dimer dynamics and its applications to saddle point search. {\it SIAM J. Numer. Anal.} 50 (2012), 1899--1921.



\bibitem{ZhaSINUM} L. Zhang, P. Zhang, X. Zheng, Error estimates for Euler discretization of high-index saddle dynamics. {\it SIAM J. Numer. Anal.} 60 (2022), 2925--2944.

\bibitem{ZhaZhe} L. Zhang, X. Zheng, S. Zhu, Numerical analysis for saddle dynamics of some semilinear elliptic problems. {\it SIAM J. Numer. Anal.} 64 (2026), 755--778.

\bibitem{ZhaLi} X. Zhang, Z. Li,  L. Zhang, Constructing a holistic map of cell fate decision by hyper solution landscape. {\it Cell Systems}, 17 (2026), p. 101562.








\end{thebibliography}
\end{document}